\documentclass[11pt]{article}

\usepackage{fullwidth}
\usepackage{xypic}
\usepackage[top=1.5in, bottom=1.5in, left=1in, right=1in]{geometry}
\usepackage{amsgen}
\usepackage{amsmath}
\usepackage{amstext}
\usepackage{amsbsy}
\usepackage{amsopn}
\usepackage{amsfonts}
\usepackage{amssymb}
\usepackage{eepic}
\usepackage{graphicx}
\usepackage{epsf}
\usepackage{pstricks}
\xyoption{all}

\def\Box{\square}
\def\edge{\relbar\joinrel\relbar}
\def\longedge{\relbar\joinrel\relbar\joinrel\relbar\joinrel\relbar}
\def\mapright#1{\smash{\mathop{\longrightarrow}\limits^{#1}}}
\def\tra#1{\smash{\mathop{\mid\kern
-1pt\joinrel\relbar\joinrel\relbar}\limits^{*}_{#1}}}
\def\longtra#1{\smash{\mathop{\mid\kern
-1pt\joinrel\relbar\joinrel\relbar\joinrel\relbar}\limits^{*}_{#1}}}
\def\vlongtra#1{\smash{\mathop{\mid\kern
-1pt\joinrel\relbar\joinrel\relbar\joinrel\relbar\joinrel\relbar}\limits^{*}_{#1}}}
\def\vvlongtra#1{\smash{\mathop{\mid\kern
-1pt\joinrel\relbar\joinrel\relbar\joinrel\relbar\joinrel\relbar\joinrel\relbar}\limits^{*}_{#1}}}
\def\vvvlongtra#1{\smash{\mathop{\mid\kern
-1pt\joinrel\relbar\joinrel\relbar\joinrel\relbar\joinrel\relbar\joinrel\relbar\joinrel\relbar}\limits^{*}_{#1}}}
\def\etra#1{\smash{\mathop{\mid\kern
-1pt\joinrel\relbar\joinrel\relbar}\limits_{#1}}}

\def\Rw{\Rightarrow}
\def\oo{\overline}
\def\wt{\widetilde}
\def\wh{\widehat}

\def\L{{\cal{L}}}

\def\N{\mathbb{N}}

\def\pure{\mbox{pure}}
\def\lk{\mbox{lk}}
\def\cl{\mbox{Cl}}
\def\im{\mbox{Im}\,}
\def\ord{\mbox{Ord}}

\def\fct{\mbox{fct}}
\def\nbh{\mbox{nbh}}

\def\dim{\mbox{dim}}

\def\ker{\mbox{Ker}\,}

\def\min{\mbox{min}}

\def\flatx{\mbox{Fl}}

\def\H{{\cal{H}}}
\def\BoR{{\cal{BR}}}

\def\Z{\mathbb{Z}}

\def\p{\varphi}

\def\inv{^{-1}}

\def\bi{\begin{itemize}}
\def\ei{\end{itemize}}
\def\beq{\begin{equation}}
\def\eeq{\end{equation}}

\def\J{{\cal{J}}}

\newtheorem{T}{Theorem}[section]
\newcommand{\bt}{\begin{T}}
\newcommand{\et}{\end{T}}
\newcommand{\ftd}{$\square$\end{T}}

\newtheorem{Proposition}[T]{Proposition}
\newcommand{\bp}{\begin{Proposition}}
\newcommand{\ep}{\end{Proposition}}
\newcommand{\fpd}{$\square$\end{Proposition}}

\newtheorem{Lemma}[T]{Lemma}
\newcommand{\bl}{\begin{Lemma}}
\newcommand{\el}{\end{Lemma}}
\newcommand{\fld}{$\square$\end{Lemma}}

\newtheorem{Corol}[T]{Corollary}
\newcommand{\bc}{\begin{Corol}}
\newcommand{\ec}{\end{Corol}}
\newcommand{\fcd}{$\square$\end{Corol}}

\newtheorem{Result}[T]{Result}
\newcommand{\br}{\begin{Result}}
\newcommand{\er}{\end{Result}}
\newcommand{\frd}{$\square$\end{Result}}

\newtheorem{Example}[T]{Example}
\newcommand{\be}{\begin{Example}}
\newcommand{\ee}{\end{Example}}

\newtheorem{Problem}[T]{Problem}
\newcommand{\bq}{\begin{Problem}}
\newcommand{\eq}{\end{Problem}}

\newcommand{\proof}
   {\par\medbreak\noindent{\bf Proof}.\enspace}

\newcommand{\qed}{
$\Box$
\par\bigbreak}

\normalbaselines
\def\abstract#1{\par\bigskip
\begingroup\small
\baselineskip=12truept
\begin{center}ABSTRACT\end{center}
\par\medskip\par\noindent
\null\hfill\hbox{\vbox{\hsize=5truein\noindent#1}}
\hfill\null\par\endgroup\par}

\title{On the topology of a boolean representable simplicial complex} 
\author{{\bf Stuart Margolis, John Rhodes and Pedro V. Silva}}

\date{\today}

\begin{document}
\maketitle

\begin{center}\small

\bigskip

\end{center}

\abstract{It is proved that fundamental groups of boolean
  representable simplicial complexes are free and the rank is
  determined by the number and nature of the connected components of
  their graph of 
  flats for dimension $\geq 2$. In the case of dimension 2, it is
  shown that boolean representable simplicial complexes have the
  homotopy type of a wedge of spheres of dimensions 1 and 2. Also in
  the case of dimension 2, necessary and sufficient 
  conditions for shellability and being sequentially Cohen-Macaulay
  are determined. Complexity bounds are provided for all the
  algorithms involved.}

\section{Introduction}

In a series of three papers \cite{IR1,IR2,IR3}, Izhakian and Rhodes
introduced the concept of boolean representation for various algebraic
and combinatorial structures. These ideas were inspired by previous
work by Izhakian and Rowen on supertropical matrices (see
e.g. \cite{Izh2,IRow2,IRow,IRow22}), and were subsequently developed by
Rhodes and Silva in a recent monograph, devoted to boolean
representable simplicial complexes \cite{RSm}.

The original approach was to consider matrix representations over the
superboolean semiring $\mathbb{SB}$, using appropriate notions of vector
independence and rank. Writing $\mathbb{N} = \{ 0,1,2,\ldots \}$, we
can define $\mathbb{SB}$ as the quotient of $(\mathbb{N},+,\cdot)$ (usual
operations) by the congruence which identifies all integers $\geq 2$.
In this context, boolean representation refers to matrices using only
0 and 1 as entries.

In this paper, we view (finite) simplicial complexes under two
perspectives, geometric and combinatorial. It is well known that each
structures determines the other (see e.g. \cite[Section A.5]{RSm}). 

As an alternative to matrices, boolean representable simplicial complexes
can be characterized by 
means of their lattice of flats. The lattice of flats plays a fundamental
role in matroid theory but is not usually considered for arbitrary
simplicial complexes, probably due to the fact that, unlike the
matroid case, the structure of a simplicial complex cannot in general be
recovered from its lattice of flats. However, this is precisely what
happens with boolean representable simplicial complexes. If $\H\; =
(V,H)$ is a simplicial complex and $\flatx\H$ denotes its lattice of
flats, then $\H$ is boolean representable if and only if $H$ equals
the set of transversals of the successive differences for chains in
$\flatx\H$. This implies in particular that all (finite) matroids are boolean
representable.

In this paper we begin the study of the topology of boolean
representable simplicial complexes (BRSC). 

As any finitely presented group can be the fundamental group of a
2-dimensional simplicial complex (see e.g. \cite[Theorem 7.45]{Rot}),
the problem of understanding the homotopy type of an arbitrary
simplicial complex is hopeless.  

However, for matroids, the topology is very restricted. Indeed, it is
known that a matroid is pure shellable \cite{Bjo}. This implies
that a matroid of rank $r$ has the homotopy type of a wedge of $r-1$
dimensional spheres, the number of which is then the rank of its
unique non-trivial homology group. This latter number has a number of
combinatorial interpretations \cite{Bjo}. In particular, a matroid of
dimension at least 2 has a trivial fundamental group.  

One of the main results of this paper is to show that the fundamental
group of a BRSC is a free group. We give a precise formula for the
rank of this group in terms of the number and nature of the connected
components of its graph 
of flats \cite{RSm}. In the simple case, this rank is equivalently a
function of the number of 
connected components of the proper part of its lattice of flats.  

For 2 dimensional BRSCs, we completely
characterize shellable complexes, showing that these are precisely the
sequentially Cohen-Macauley complexes
\cite{BWW}. Although not every 2 dimensional BRSC is shellable, we
prove that every 2 dimensional BRSC has the homotopy type of a wedge
of 1-spheres and 2-spheres.

We consider the connection to EL-labelings \cite{Bjo}
of the lattice of flats
and give an example of a shellable 2-dimensional complex whose lattice
of flats is not EL-labelable. 



The paper is organized as follows. In Section 2 we present 
basic notions and results needed in the paper. In Section 3 we show
that the fundamental group of a boolean representable simplicial
complex is always free, and provide an exact formula to compute its
rank for dimension $\geq 2$,
using the  graph of flats. We also prove that any 2 dimensional BRSC
has the homotopy type of a wedge of 1-spheres and 2-spheres.

For higher degree homotopy groups, the situation is of course much
harder, and we limit the discussion to shellability in dimension 2. We
note that in \cite{RSm} we had characterized shellability for simple
boolean representable complexes of dimension 2. We are now able to
deal with the non simple case, and to assist us on this reduction we
use the concept of simplification in Section 4. Then Section 5 
is devoted to characterizing shellability for boolean representable
simplicial complexes of dimension 2. For such complexes, it is also shown
that the shellable complexes are precisely the sequentially
Cohen-Macaulay complexes.

In Section 6, we consider the concept of the order complex of a
lattice $L$. The vertices of the order complex are the elements of the
proper part of $L$, i.e. $L^* = L \setminus \{ 0,1\}$,
and its faces are the chains of $L^*$.
We show that, given a boolean representable simplicial
complex $\H$, if the order complex of ${\rm Fl}\H$ is shellable, so is
$\H$. The converse turns out to be false.

In the matroid case, (some) shellings can be obtained from EL-labelings of
the lattice of flats (which is always geometric and thus has an
EL-labeling by a theorem of Bj\"orner \cite{Bjo2}). We show that, for
arbitrary shellable pure boolean representable
simplicial complexes of dimension 2, the lattice of flats does not
necessarily admit an EL-labeling.

Finally, Section 7 discusses the complexity of several algorithms
designed to compute fundamental groups, decide shellability (for
dimension 2) and compute shellings and Betti numbers. Although the
number of potential flats in a simplicial complex with $n$ vertices is
$2^n$ and therefore exponential, we achieve
polynomial bounds for all algorithms when the dimension of the
simplicial complexes is fixed.

\section{Preliminaries}

All lattices and simplicial complexes in this paper are assumed to be
finite. Given a set $V$ and $n \geq 0$, we denote by $P_n(V)$
(respectively $P_{\leq n}(V)$) the set
of all subsets of $V$ with precisely (respectively at most) $n$
elements. The {\em kernel} of a mapping $\p:V \to W$ is the relation
$$\ker\p = \{ (a,b) \in V \times V \mid a\p = b\p\}.$$ 

A (finite) simplicial complex is a structure of the form $\H\; =
(V,H)$, where $V$ is a finite nonempty set and $H \subseteq 2^V$
contains $P_1(V)$ and is closed under taking subsets. The elements of
$V$ and $H$ are called respectively {\em vertices} and {\em faces}. To
simplify notation, we shall often denote a face $\{
x_1,x_2,\ldots,x_n\}$ by $x_1x_2\ldots x_n$. 

A face of $\H$ which is maximal with respect to
inclusion is called a {\em facet}. We denote by ${\rm fct}\H$ the set
of facets of $\H$.

The {\em dimension} of a face $I \in H$ is $|I|-1$. An $i$-{\em face}
(respectively $i$-{\em facet}) is a face (respectively facet) of
dimension $i$. We may refer to 0-faces and 1-faces as vertices and edges.

We say that $\H$ is:
\bi
\item {\em simple} if $P_2(V) \subseteq H$;
\item {\em pure} if all the facets
of $\H$ have the same dimension. 
\ei
The dimension of $\H$, denoted by $\dim\H$, is the maximum dimension of a
face(t) of $\H$.



Given  $Q \in
H \setminus \{ V \}$, we define the {\em link} $\lk(Q)$ to be the  
simplicial complex $(V/Q,H/Q)$, where
$$H/Q = \{ X \subseteq V\setminus Q \mid X \cup Q \in
H\}\quad\mbox{and}\quad
V/Q = \bigcup_{X \in H/Q} 2^X.$$
Here it is convenient to admit a simplicial complex to have an empty set
of vertices.

A simplicial complex $\H\; = (V,H)$ is called a {\em matroid} if it
satisfies the {\em exchange property}:
\bi
\item[(EP)]
For all $I,J \in H$ with $|I| = |J|+1$, there exists some
  $i \in I\setminus J$ such that $J \cup \{ i \} \in H$.
\ei

A simplicial complex $\H\; = (V,H)$ is {\em
  shellable} if we can order its facets as $B_1, \ldots, B_t$ so
that, for $k = 2, \ldots,t$, the following condition is satisfied: if
$I(B_k) = (\cup_{i=1}^{k-1} 2^{B_i}) \cap 2^{B_k}$, then
$$(B_k,I(B_k))\mbox{ is pure of dimension }|B_k|-2$$
whenever $|B_k| \geq 2$.
Such an ordering is
called a {\em shelling}. In the literature, this is called non-pure
shellability and was first defined by Bj\"orner and Wachs \cite{BW,BW2}. 

Given an $R \times V$ matrix $M$ and $Y \subseteq R$, $X \subseteq V$,
we denote by $M[Y,X]$ the submatrix of $M$ obtained by deleting all
rows (respectively columns) of $M$ which are not in $Y$ (respectively $X$).

A boolean matrix $M$ is {\em lower unitriangular} if it is of the form
$$\left(
\begin{matrix}
1&&0&&0&&\ldots&&0\\
?&&1&&0&&\ldots&&0\\
?&&?&&1&&\ldots&&0\\
\vdots&&\vdots&&\vdots&&\ddots&&\vdots\\
?&&?&&?&&\ldots&&1
\end{matrix}
\right)
$$

Two matrices are {\em congruent} if we can transform one into the other
by independently permuting rows/columns. A boolean matrix is {\em nonsingular}
if it is congruent to a lower unitriangular matrix.

Given an $R \times V$ boolean matrix $M$, we say that the
subset of columns $X \subseteq V$ is $M$-{\em independent} if there exists
some $Y \subseteq R$ such that $M[Y,X]$ is nonsingular.

A simplicial complex $\H\; = (V,H)$ is {\em boolean representable} if
there exists some boolean matrix $M$ such that
$H$ is the set of all $M$-independent subsets of $V$.

We denote by $\BoR$
the class of all (finite) boolean representable simplicial
complexes. All matroids are boolean representable \cite[Theorem
  5.2.10]{RSm}, but the converse is not true.

We say that $X
\subseteq V$ is a {\em flat} of $\H$ if
$$\forall I \in H \cap 2^X \hspace{.2cm} \forall p \in V \setminus X
\hspace{.5cm} I \cup \{ p \} \in H.$$
The set of all flats of $\H$ is denoted by 
$\flatx\H$. Note that $V, \emptyset \in \flatx\H$ in all cases.

Clearly, the intersection of any set of flats (including $V =
\cap\emptyset$) is still a flat. If we order $\flatx\H$ by inclusion,
it is then a $\wedge$-semilattice. Since $\flatx\H$ is finite, it
follows that it is indeed a lattice (with the determined join), the
{\em lattice of flats} of $\H$. 

We say that $X$ is a {\em transversal of the
successive differences} for a chain of subsets
$$A_0 \subset A_1 \subset \ldots \subset A_k$$
if $X$ admits an enumeration $x_1,\ldots , x_k$ such that $x_i \in A_i
\setminus A_{i-1}$ for $i = 1,\ldots,k$. 

Let $\H\; = (V,H)$ be a simplicial complex. If $X \subseteq V$ is a
{\em transversal of the successive differences} for a chain 
$$F_0 \subset F_1 \subset \ldots \subset F_k$$
in $\flatx \H$, it follows easily by induction that $x_1x_2
\ldots x_i \in H$ for $i = 0,\ldots,k$. In particular, $X \in H$.

It follows from \cite[Corollary 5.2.7]{RSm} that $\H$ is boolean
representable if and only if every $X 
\in H$ is a transversal of the 
successive differences for a chain in $\flatx \H$. 

The lattice $\flatx\H$ induces a closure operator on $2^V$ defined by
$$\oo{X} = \cap\{ F \in \flatx \H \mid X \subseteq F \}$$
for every $X \subseteq V$. 

By \cite[Corollary 5.2.7]{RSm}, $\H\; = (V,H)$ is boolean
representable if and only if every $X \in H$ admits an enumeration
$x_1,\ldots, x_k$ satisfying
\beq
\label{derby3}
\oo{x_1} \subset \oo{x_1x_2} \subset\ldots \subset \oo{x_1\ldots x_k}.
\eeq
Thus, given $p,q \in V$ distinct, we have 
\beq
\label{cheta}
pq \notin H \mbox{ if and only } \oo{p} = \oo{pq} = \oo{q}.
\eeq
This fact will be often used throughout the text with no explicit reference.
From (\ref{cheta}) we can deduce that
\beq
\label{cheta1}
\oo{p} = \{ q \in V \mid \oo{q} = \oo{p} \}.
\eeq
Indeed, let $F = \{ q \in V \mid \oo{q} = \oo{p} \}$. Since $p \in F
\subseteq \oo{p}$, it suffices to show that $F \in \flatx\H$. Let $I
\in H \cap 2^F$ and $a \in V \setminus F$. In view of (\ref{cheta}),
we may assume that $I = \{ q \}$. Since $\oo{a} \neq \oo{q}$, we get
$qa \in H$ also by (\ref{cheta}). Thus $F \in \flatx\H$ and
(\ref{cheta1}) holds.

Let $\J\; = (V, J)$ be a simplicial complex. We recall the definitions
of the (reduced) homology groups of $\J$ (see e.g. \cite{Hat}).

If $\J$ has $s$ connected components, it is well known that the {\em
  0th homology group} $H_0(\J)$ is isomorphic to 
the free abelian group of rank $s$. For dimension $k \geq 1$, we
proceed as follows.

Fix a total ordering of $V$.
Let $C_k(\J)$ denote the free abelian group on $J \cap P_{k+1}(V)$,
that is, all the formal sums of the 
form $\sum_{i \in I} 
n_iX_i$ with $n_i \in \Z$ and $X_i \in J \cap P_{k+1}(V)$
(distinct). Given $X \in J \cap P_{k+1}(V)$, write $X = x_0x_1\ldots x_k$
with $x_0 < \ldots < x_k$. We define 
$$X\partial_k = \sum_{i= 0}^k (-1)^i (X\setminus\{ x_i\}) \in
C_{k-1}(\J)$$
and extend this by linearity to a homomorphism $\partial_k:C_{k}(\J)
\to C_{k-1}(\J)$ (the $k$th {\em boundary map} of $\J$). Then the {\em
  $k$th homology group} of $\J$ is defined as the quotient
$$H_k(\J) = \ker\partial_k / \im\partial_{k+1}.$$

The {\em $0$th reduced homology group} of $\J$, denoted by
$\tilde{H}_0(\J)$, is isomorphic to the free abelian group of rank $s-1$, where
$s$ denotes the number of connected components of $\J$.
For $k \geq 1$, the {\em $k$th reduced homology group} of $\J$,
denoted by $\tilde{H}_k(\J)$ coincides with the $k$th homology group.

A {\em wedge} of spheres $S_1,\ldots, S_m$ (of possibly different
dimensions) is a 
topological space obtained by identifying $m$ points $s_i \in S_i$ for
$i = 1,\ldots,m$.  

Given a group $G$ and $X \subseteq G$, we denote by $\langle X
\rangle$  (respectively $\langle\langle X
\rangle\rangle$) the subgroup (respectively normal subgroup) of $G$
generated by $X$. 

We denote by $F_A$ the free group on an alphabet $A$.
A {\em group presentation} is a formal expression of
the form $\langle A \mid R\rangle$, where $A$ is an alphabet and $R
\subseteq F_A$. It defines the group $F_A/\langle\langle
R\rangle\rangle$,
and is said to be a presentation for any group
isomorphic to this quotient.  

Given a (finite) alphabet $A$, we denote by $A^+$ the {\em free
  semigroup} on $A$ (finite nonempty words on $A$, under
concatenation). Given a partial order on $A$, we define the {\em
  lexicographic order} on $A^+$ as follows. Given $a_1,
\ldots,a_k,a'_1,\ldots,a'_m \in A$, we write $a_1 \ldots a_k < a'_1
\ldots a'_m$ if one of the following conditions holds:
\bi
\item
$k < m$ and $a_i = a'_i$ for $i = 1,\ldots,k$;
\item
there exists some $i \leq \min \{ k,m \}$ such that $a_1 = a'_1,\;
\ldots, \; a_{i-1}
= a'_{i-1}, \; a_i < a'_i$.
\ei

\section{The fundamental group}
\label{fugr}

Let $\H\; = (V,H)$ be a simplicial complex. The {\em graph} of
$\H$ is the truncation $(V,H \cap P_{\leq 2}(V))$. We say that
$\H$ is {\em connected} if its graph is connected. We say that $T \subseteq H 
\cap P_2(V)$ is a spanning tree of $\H$ if it is a spanning tree of its graph. 

\bl
\label{connected}
Let $\H\; = (V,H)$ be a boolean representable simplicial complex. Then
$\H$ is connected unless $H = P_1(V)$ and $|V| > 1$.
\el

\proof
Obviously, $\H$ is disconnected if $H = P_1(V)$ and $|V| > 1$, and
connected if $|V| = 1$. Hence
we may assume that $pq \in H$ for some distinct $p,q \in V$. 

Let $M$ be an $R \times V$ boolean matrix representing $\H$. It
follows from $pq \in H$ that $M[R,p] \neq M[R,q]$. Thus, for every $v
\in V$, we have either $M[R,v] \neq M[R,p]$ or $M[R,v] \neq M[R,q]$,
implying that $vp$ or $vq$ is an edge of $H$.
Therefore $\H$ is connected.
\qed

Note that, if we consider the {\em geodesic distance} on the graph of
a boolean representable simplicial complex of
dimension $\geq 2$ (the distance 
between two vertices is the length of the shortest path connecting
them), it follows from the above proof that the distance between any
two vertices is at most 2.

It is well known that the geometric realization $||\H||$ of a simplicial
complex, a subspace of some euclidean space $\mathbb{R}^n$, is unique
up to homeomorphism. For details, see e.g. \cite[Appendix A.5]{RSm}. 

Given a point $v_0 \in ||\H||$, the {\em
  fundamental group} $\pi_1(||\H||,v_0)$ 
is the group having as
elements the homotopy equivalence classes of closed paths
$$
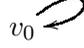
\xymatrix{
v_0 \ar@(ur,r)
}$$
the product being determined by the concatenation of paths. 

If $\H$ is connected, then $\pi_1(||\H||,v_0) \cong \pi_1(||\H||,w_0)$
for all points $v_0,w_0$ in $||\H||$, hence we may use the notation
$\pi_1(||\H||)$ without ambiguity. We produce now a presentation for
$\pi_1(||\H||)$. This combinatorial description is also known as the
{\em edge-path group} of $\H$ (for details on the
fundamental group of a simplicial complex, see \cite{Spa}). 

We fix a spanning tree $T$ of $\H$ and we define 
$$A = \{ a_{pq} \mid pq \in H \cap P_2(V)\},$$
$$R_T = \{ a_{qp}a_{pq}\inv \mid pq \in H \cap P_2(V) \} \cup \{
a_{pq}a_{qr}a_{pr}\inv \mid pqr \in H \cap P_3(V) \} 
\cup \{ a_{pq} \mid pq \in T \}.$$
From now on, we view $\pi_1(||\H||)$ as the group
defined by the group presentation 
\beq
\label{pre}
\langle A \mid R_T \rangle.
\eeq
We denote by $\theta:F_A \to \pi_1(||\H||)$ the canonical homomorphism.
We note that the six relators induced by a single 2-face $pqr$
(corresponding to different enumerations of the vertices) are all
equivalent to $a_{pq}a_{qr}a_{pr}$: each one of them is a conjugate of
either $a_{pq}a_{qr}a_{pr}$ or its inverse. 

Given a boolean representable connected simplicial
complex $\H\; = (V,H)$, the {\em graph of flats} $\Gamma\flatx\H$ has
vertex set $V$ and edges $p \edge q$ whenever $p \neq q$ and $\oo{pq}
\subset V$. 


\bl
\label{count}
Let $\H\; = (V,H)$ be a boolean representable connected 
simplicial complex. Let $u,v \in V$ belong to distinct connected
components of $\Gamma{\rm Fl}(\H)$. Then $uv \in H$.
\el

\proof
Since $|V| > 1$ and $\H$ is connected, there exists some $pq \in H
\cap P_2(V)$. Suppose that $uv \notin H$. By (\ref{cheta}), we get
$\oo{u} = \oo{uv} = \oo{v}$. Since there is no edge $u \edge v$ in
$\Gamma\flatx\H$, we get $\oo{u} = V$. By (\ref{cheta1}), we
get $\oo{p} = \oo{q} = \oo{u} = V$. In view of (\ref{cheta}), this
contradicts $pq \in H$.
%
%
%
%
\qed




Let $C$ be a connected component of $\Gamma{\rm Fl}(\H)$. If $H \cap
P_2(C) \neq \emptyset$, we shall say that $C$ is 
$H$-{\em nontrivial}. Otherwise, we say that $C$ is $H$-{\em trivial}.
The {\em size} of $C$ is its number of vertices.

If $\H$ is a connected simplicial complex of dimension $\leq
1$ (i.e. a graph), then (\ref{pre}) is a presentation of a free group,
its rank equal to the number of edges of the graph that are not in $T$. 

The next result shows that the graph of flats and the size of its
$H$-disconnected components determines completely
the fundamental group for dimension $\geq 2$.

\bt
\label{fung}
Let $\H$ be a boolean representable simplicial
complex of dimension $\geq 2$. Assume that $\Gamma{\rm Fl}\H$ has $s$
$H$-nontrivial connected components and $r$ $H$-trivial connected
components of sizes $f_1,\ldots,f_r$. 
Then $\pi_1(||\H||)$ is a free group of rank 
$$\binom{s+f_1+\ldots +f_r-1}{2} - \sum_{i=1}^r \binom{f_i}{2},$$
or equivalently, 
$$\binom{s-1}{2} +(s-1)(f_1+\ldots +f_r) + \sum_{1 \leq i < j \leq r} f_if_j.$$
\et

\proof
Let $\H\; = (V,H)$ and $\Gamma = \Gamma\flatx\H$. 
Since $\H$ has dimension $\geq 2$, there
exists some $xyz \in H \cap P_3(V)$. Since $\H$ is boolean
representable, we may assume by (\ref{derby3}) that $\oo{yz} \subset
V$, hence $y \edge z$ is an edge of $\Gamma$.
In view of (\ref{cheta}), we may also assume that $y \notin \oo{z}$.

Let
$$Z = \{ p \in V\setminus \{ z \} \mid pz \in H\}.$$
Note that $y \in Z$. Now let
$$T = \{ pz \mid p \in Z\} \cup \{ yq \mid q \in V\setminus (Z \cup \{
z \})\}.$$
We claim that $T$ is a spanning tree of $\H$.

Indeed, suppose that $q \in V\setminus (Z \cup \{
z \})$. Then $qz \notin H$ and so $\oo{q} = \oo{qz} = \oo{z}$.
Since $y \notin \oo{z}$, we get $y \notin \oo{q}$, hence $yq \in H$
and so $T \subseteq H \cap 
P_2(V)$. Now $T$ has precisely $|V|-1$ edges and every vertex of $V$
occurs in some edge of $T$. Therefore $T$ is a spanning tree of $\H$.

We consider now the finite presentation (\ref{pre}) of $\pi_1(||\H||)$
induced by the spanning tree $T$. Our goal is to use a sequence of
{\em Tietze transformations} (see \cite{LS}) to obtain a presentation
that can be seen to be that of the free group in the statement of the
theorem. This requires some preliminary work.

Let $\theta:F_A \to \pi_1(||\H||)$
denote the canonical homomorphism. We show that
\beq
\label{fung2}
pq \in E(\Gamma) \cap H \Rw a_{pq}\theta = 1.
\eeq

Suppose first that $z \notin \oo{pq}$. Then $pqz \in H$, hence $p,q
\in Z$ and we get
$$a_{pq}\theta = (a_{zp}a_{pq}a_{zq}\inv)\theta = 1.$$
Thus we may assume that $z \in \oo{pq}$.

Suppose that $y \notin \oo{pq}$. Then $pqy \in H$. We claim that 
\beq
\label{fung3}
a_{yp}\theta = a_{yq}\theta = 1.
\eeq

If $p \in V \setminus Z$, then $yp \in T$ and so $a_{yp}\theta = 1$.
If $p \in Z$, then $pz \in T$. Since $\oo{pz}
\subseteq \oo{pq}$ yields $y \notin \oo{pz}$, we get $yzp \in H$ and
so
$$a_{yp}\theta = (a_{yz}a_{zp})\theta = 1.$$

Similarly, $a_{yq}\theta = 1$ and so (\ref{fung3}) holds.

Now $pqy \in H$ yields
$$a_{pq}\theta = (a_{py}a_{yq})\theta = (a_{yp}\inv a_{yq})\theta = 1.$$

So finally we may assume that $z,y \in \oo{pq}$. Let $v \in V\setminus
\oo{pq}$. We prove that $a_{pv}\theta = 1$ by considering two cases.
If $\oo{p} \neq \oo{z}$, then $pzv \in H$ and so
$a_{pv}\theta = (a_{pz}a_{zv})\theta = 1$. Hence we assume that
$\oo{p} = \oo{z}$. Now $yzv \in H$ yields $a_{yv}\theta =
(a_{yz}a_{zv})\theta = 1$, and $pyv \in H$ (which holds since $\oo{p}
= \oo{z}$ implies $\oo{p} \neq \oo{y}$) yields $a_{pv}\theta =
(a_{py}a_{yv})\theta = 1$ (since $py \in T$). 

Hence $a_{pv}\theta = 1$ and by symmetry also $a_{qv}\theta =
1$. Finally, $pqv \in H$ yields $a_{pv}\theta =
(a_{pq}a_{qv})\theta$ and thus $a_{pq}\theta = 1$. 
Therefore (\ref{fung2}) holds.









Let $C_1, \ldots,C_s$ (respectively $C'_1,\ldots,C'_r$) denote the 
$H$-nontrivial (respectively $H$-trivial)
connected components of $\Gamma$. We assume also that $C'_i$ has size
$f_i$ for $i = 1,\ldots,r$.

We say that two vertices $p,q \in C_i$ are $H$-{\em connected} if
there exists a path
$$p = p_0 \edge p_1 \edge \ldots \edge p_n = q$$
in $C_i$ with $n \geq 0$ and $p_{j-1}p_j \in H$ for $j = 1,\ldots,n$.

We claim that
\beq
\label{fung6}
pq \in H \cap P_2(C_i) \Rw p\mbox{ and $q$ are $H$-connected}
\eeq
holds for $i = 1,\ldots,s$.

Let $d$ denote the geodesic distance on $C_i$.
We show that $p,q \in C_i$ are $H$-connected using induction on $d(p,q)$.

The case $d(p,q) \leq 1$ is trivial, hence we assume that
$d(p,q) = n > 1$ and (\ref{fung5}) holds for closer vertices. Take $p',p''
\in C_i$ such that $d(p,p') = n-2$ and $d(p',p'') = d(p'',q) = 1$:
$$p \longedge p' \edge p'' \edge q$$ 

Suppose that $p''q \notin H$. Then $\oo{p''} = \oo{p''q} = \oo{q}$. It
follows that $\oo{p'q} = 
\oo{p'p''} \subset V$ and so there exists an edge $p' \edge q$ in
$\Gamma$, contradicting $d(p,q) = n$.

Thus $p''q \in H$. Since $d(p,q) > 1$, we have $p \notin
\oo{p''q} \subset V$. Hence $pp'' \in H$. But $d(p,p'') = n-1$, so by
the induction hypothesis $p$ and $p''$ are $H$-connected. Since $p''q
\in H$, it follows that $p$ and $q$ are $H$-connected. Therefore
(\ref{fung6}) holds.

We show next that
\beq
\label{fung5}
pq \in H \cap P_2(C_i) \Rw a_{pq}\theta = 1
\eeq
holds for $i = 1,\ldots,s$.

We use induction on $d(p,q)$.
The case $d(p,q) = 1$ follows from (\ref{fung2}), hence we assume that
$d(p,q) = n > 1$ and (\ref{fung5}) holds for closer vertices. Take $p',p''
\in C_i$ as in the proof of (\ref{fung6}). By that same proof, we must have
$p''q \in H$. Since $d(p,q) > 1$, we have $p \notin
\oo{p''q}$. Hence $pp''q \in H$ and so $pp'', p''q \in H$. By the induction
hypothesis, we get $a_{pp''}\theta = a_{p''q}\theta = 1$.
But now $pp''q \in H$ yields $a_{pq}\theta = (a_{pp''}a_{p''q})\theta =
1$. Therefore
(\ref{fung5}) holds.

Now we may use (\ref{fung5}) to simplify the group presentation
$\langle A \mid R_T\rangle$. In view of (\ref{fung5}), we start by
adding as relators all the 
$a_{pq} \in A$ such 
that $p,q$ belong to the same $C_i$. 

For $i = 1,\ldots,s$, we fix some vertex $c_i \in C_i$. We may assume
without loss of generality that $c_1 = z$. Given $p \in V$, we write
$\wh{p} = c_i$ if $p \in C_i$. We define
$$\begin{array}{lll}
R'&=&\{ a_{qp}a_{pq}\inv \mid pq \in H \cap P_2(V) \} \cup 
\{ a_{pq} \mid pq \in T \}\\
&\cup&\{ a_{pq} \mid pq \in H \cap P_2(C_i),\; i \in \{ 1,\ldots,s \} \}\\ 
&\cup&\{ a_{pq}a_{\wh{p}\wh{q}}\inv \mid p \in C_i,\: q \in C_j,\; i,j
\in \{ 1,\ldots,s \}, \; i \neq j \}\\
&\cup&\{ a_{pq}a_{\wh{p}q}\inv \mid p \in C_i,\: q \in C'_j,\; i
\in \{ 1,\ldots,s \}, \; j \in \{ 1,\ldots, r \} \}\\
&\cup&\{ a_{pq}a_{p\wh{q}}\inv \mid p \in C'_j,\: q \in C_i,\; i
\in \{ 1,\ldots,s \}, \; j \in \{ 1,\ldots, r \} \}.
\end{array}$$
In view of Lemma \ref{count}, $R'$ is well defined.
We show that $\langle\langle R' \rangle\rangle = \langle\langle R_T
\rangle\rangle$. 

We show first that $R' \subseteq \langle\langle R_T
\rangle\rangle$. In view of (\ref{fung5}), we only need to discuss the
last three terms of the union.

We start by proving that
\beq
\label{fung8}
a_{pq}\theta = a_{\wh{p}q}\theta
\eeq
whenever $p \in C_i$ and $q \notin C_i$.
We may assume that $p \neq \wh{p}$. By (\ref{fung6}),
there exists a path
$$p = p_0 \edge p_1 \edge \ldots \edge p_n = \wh{p}$$
in $C_i$ with $n \geq 1$ and $p_{k-1}p_k \in H$ for $k = 1,\ldots,n$.
By Lemma \ref{count}, we have $p_kq \in H$ for every $k$. Also
$\oo{p_{k-1}p_k} \subset V$ for $k = 1,\ldots,n$. Suppose that $q
\in \oo{p_{k-1}p_k}$. Then $\oo{p_kq} \subset V$ and $q \in C_i$, a
contradiction. Hence $q \notin \oo{p_{k-1}p_k}$. Since $p_{k-1}p_k \in
H$, it follows that $p_{k-1}p_kq \in
H$ and in view of (\ref{fung5}) we get
$$a_{p_{k-1}q}\theta = (a_{p_{k-1}p_k}a_{p_{k}q})\theta =
a_{p_{k}q}\theta.$$
Now (\ref{fung8}) follows by transitivity.

Similarly,
\beq
\label{fung13}
a_{pq}\theta = a_{p\wh{q}}\theta
\eeq
whenever $q \in C_i$ and $p \notin C_i$.

Finally, if $p \in C_i$ and $q \in C_j \neq C_i$, we may apply
(\ref{fung8}) and (\ref{fung13}) to get
$a_{pq}\theta = a_{\wh{p}q}\theta = a_{\wh{p}\wh{q}}\theta$.
Therefore $R' \subseteq \langle\langle R_T
\rangle\rangle$ and so $\langle\langle R'
\rangle\rangle \subseteq \langle\langle R_T
\rangle\rangle$.

To prove the opposite inclusion, let $\theta':F_A \to F_A/\langle\langle R'
\rangle\rangle$ denote the canonical homomorphism. It suffices to show
that $(a_{pq}a_{qr}a_{pr}\inv)\theta' = 1$ for every $pqr \in H \cap
P_3(V)$.

Since $\H$ is boolean representable and $pqr \in H$, one of the three
elements $p,q,r$ is not in the closure of the other two. We remarked
before that each one of the six relators of $R_T$ arising from
distinct enumerations of the elements of $p,q,r$ is a conjugate of
$a_{pq}a_{qr}a_{pr}\inv$ or its inverse, hence we may assume that $r
\notin \oo{pq}$. Hence there exists an edge $p \edge q$ in $\Gamma$
and so $p,q \in C_i$ for some $i \in \{ 1,\ldots,s\}$. 

Suppose that $r \in C_i$. Since $pq,qr,pr \in H$, we get
$a_{pq}\theta' = a_{qr}\theta' = a_{pr}\theta' = 1$ and so
$(a_{pq}a_{qr}a_{pr}\inv)\theta' = 1$.

Thus we may assume that $r \notin C_i$. If $r \notin C'_1 \cup \ldots
\cup C'_r$, then
$$a_{qr}\theta' = a_{\wh{q}\wh{r}}\theta' = a_{\wh{p}\wh{r}}\theta' = 
a_{pr}\theta'.$$
The case $r \in C'_1 \cup \ldots
\cup C'_r$ is analogous.
Since $pq \in H \cap P_2(C_i)$ yields $a_{pq}\theta' = 1$, we get
$(a_{pq}a_{qr}a_{pr}\inv)\theta' = 1$. Therefore $\langle\langle R'
\rangle\rangle = \langle\langle R_T
\rangle\rangle$. 

Now we simplify the presentation $\langle A \mid R' \rangle$ by means
of further Tietze transformations.

The third term of the union in $R'$ ensures that we may omit all
generators with both indices in the same connected components, and the
three last terms allow us to restrict ourselves to generators with
indices in $\{ c_1,\ldots,c_s\} \cup C'_1 \cup \ldots \cup
C'_r$. Since $y,z \in C_1$, the second term allows us to eliminate all
the generators where $c_1 = z$ appears as index, and we may now use
the first term relators to remove half of the remaining generators,
ending up with the free group on the set
$$\begin{array}{lll}
B&=&\{ a_{c_ic_j} \mid 2 \leq i < j \leq s \}\\
&\cup&\{ a_{c_iq} \mid 2 \leq i \leq s,\; q \in C'_1\cup \ldots \cup
C'_r \}\\
&\cup&\{ a_{pq} \mid p \in C'_i,\; q \in C'_j \; 1 \leq i < j\leq r
\}
\end{array}$$

Now 
$$|B| = \binom{s-1}{2} +(s-1)(f_1+\ldots +f_r) + \sum_{1 \leq i < j \leq r} f_if_j.$$
On the other hand, we have
$$\begin{array}{lll}
\binom{s+f_1+\ldots +f_r-1}{2}&=&\frac{(s-1+f_1+\ldots
  +f_r)(s-2+f_1+\ldots +f_r)}{2}\\ 
&=&\frac{(s-1)(s-2)}{2} + (s-1)(f_1+\ldots +f_r) + \frac{(f_1+\ldots
  +f_r)(f_1+\ldots +f_r-1)}{2}\\
&=&\binom{s-1}{2} + (s-1)(f_1+\ldots +f_r) + \sum_{1 \leq i < j \leq
  r} f_if_j + \frac{\sum_{i=1}^r (f_i^2-f_i)}{2}\\
&=&\binom{s-1}{2} + (s-1)(f_1+\ldots +f_r) + \sum_{1 \leq i < j \leq
  r} f_if_j + \sum_{i=1}^r \binom{f_i}{2},
\end{array}$$
proving the theorem.
\qed

Given a lattice $L$ with top element 1 and bottom element 0, write
$L^* = L \setminus \{ 0,1 \}$ (the {\em proper part} of $L$) and
define a graph $\Delta L^* = (L^*, UH_{L^*})$, where $UH_{L^*}$
denotes the set of undirected edges in the Hasse diagram of 
$L^*$. More formally, we can define $UH_{L^*}$ as the set of all edges
$a \edge b$ such that $a$ covers $b$ in $L^*$ (i.e. $a
> b$ and there exists no $c \in L^*$ such that $a > c > b$). 

\bc
\label{simfun}
Let $\H$ be a boolean representable simple simplicial
complex of dimension $\geq 2$. 
Then $\pi_1(||\H||)$ is a free group of rank $\binom{t-1}{2}$, where
$t$ denotes the number of connected components of $\Gamma{\rm
  Fl}\H$. This number is also equal to the number of connected
components of $\Delta({\rm Fl}\H)^*$.  
\ec

\proof
If $\H\; = (V,H)$ is simple, then each $H$-trivial connected component of
$\Gamma\flatx\H$ has precisely one vertex. Hence, by Theorem
\ref{fung}, $\pi_1(||\H||)$ is a free group of rank
$\binom{t-1}{2}$.

Note that, since $\H$ is simple, then $P_1(V) \subseteq \flatx\H$ (so
all points of $\H$ belong to $(\flatx\H)^*$).

Let $p,q \in V$ be adjacent in $\Gamma\flatx\H$. Then $\oo{pq} \subset
V$ and so $\oo{pq}$ is the join of $p$ and $q$ in $\Delta({\rm
  Fl}\H)^*$. It follows that each connected component of $\Gamma{\rm
  Fl}\H$ is contained in the union of the points of some connected
component of $\Delta({\rm Fl}\H)^*$.  

On the other hand, if $F \edge F'$ is an edge of $\Delta({\rm
  Fl}\H)^*$ (say, with $F \subset F'$), then $F'$ is a clique of
$\Gamma\flatx\H$ (i.e. induces a complete subgraph). It follows easily
that the union of the points of a connected
component of $\Delta({\rm Fl}\H)^*$ belong to the same connected
component of $\Gamma\flatx\H$.

Since every connected component of $\Delta({\rm Fl}\H)^*$ contains
necessarily a point, the number of connected components must coincide
in both graphs.
\qed

We show next that free groups of rank $\binom{n}{2}$ $(n \geq 2)$
occur effectively 
as fundamental groups of boolean representable simplicial
complexes of dimension 2, even in the simple case. 

\be
\label{occur}
Let $t \geq 3$. Let $\H\; = (V,H)$ be defined by $V = \{
a_1,b_1,a_2,b_2, \ldots, a_t,b_t\}$ and
$$H = P_{\leq 2}(V) \cup \{ X \in P_3(V) \mid a_ib_i \subset X\mbox{
  for some }i \in \{ 1,\ldots,t\} \}.$$
Then $\H$ is a boolean representable simple simplicial
complex of dimension $2$ and $\pi_1(||\H||) \cong F_{\binom{t-1}{2}}$.
\ee

Indeed, it is easy to check that
$$\flatx\H = P_{\leq 1}(V) \cup \{ a_1b_1, a_2b_2, \ldots, a_tb_t,
V\},$$
hence every face of $\H$ is a transversal of the successive differences
for some chain in $\flatx\H$. Thus $\H$ is boolean
representable. Clearly, the graph of flats of $\H$ is
$$a_1 \edge b_1, \quad a_2 \edge b_2, \quad \ldots \quad a_t \edge b_t,$$
hence it possesses $t$ connected components. Therefore $\pi_1(||\H||)
\cong F_{\binom{t-1}{2}}$ by Corollary \ref{simfun}. Note also that
$\Delta(\flatx\H)^*$ is
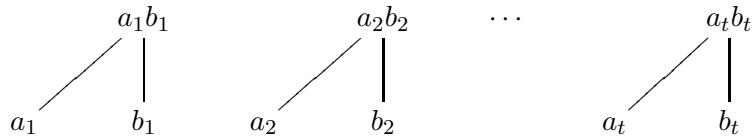
$$\xymatrix{
&{a_1b_1} \ar@{-}[dl] \ar@{-}[d] &&{a_2b_2} \ar@{-}[dl] \ar@{-}[d]
  &\cdots&&{a_tb_t} \ar@{-}[dl] \ar@{-}[d]\\
a_1&b_1&a_2&b_2&&a_t&b_t
}$$

\bigskip

By shellability of matroids, every matroid $\H\; = (V,H)$ of dimension
$d \geq 2$ has 
the homotopy type of a wedge of spheres of dimension $d$. In
particular, its fundamental group is trivial. We note that this fact
also follows from the preceding theorem, since
$\Gamma\flatx\H$ is a complete graph. Indeed, given $p,q \in V$
distinct, it is well known (see e.g. \cite[Proposition 4.2.5(ii)]{RSm}) that
$$\oo{pq} = pq \cup \{ r \in V \setminus pq \mid I \cup \{r\} \notin
H \mbox{ for some } I \in H \cap 2^{pq} \}.$$
Since every matroid is pure and $\dim\H \geq 2$, $pq$ cannot be a
facet and so $\oo{pq} \subset V$. Thus $\Gamma\flatx\H$ has a single
connected component and so $\pi_1(||\H||)$ is trivial by Theorem \ref{fung}.

Theorem \ref{fung} also yields the following consequence, one of
the main theorems of the paper.

\bt
\label{wed}
Let $\H$ be a boolean representable simplicial complex of dimension $2$. Then:
\bi
\item[(i)] the homology groups of $\H$ are free abelian;
\item[(ii)] $\H$ has the homotopy type of a wedge of 1-spheres and
  2-spheres.
\ei
\et

\proof
(i) It follows from Lemma \ref{connected} that $\H$ is connected.
By Hurewicz Theorem (see \cite{Hat}), the 1st homology group of
$\H$ is the abelianization of $\pi_1(||\H||)$, and therefore, in view
of Theorem \ref{fung}, a free abelian group of known rank. The second
homology group of any 2-dimensional simplicial complex is
$\ker\partial_2 \leq C_2(\H)$, that is, a subgroup of a free abelian
group. Therefore $H_2(\H)$ is itself free abelian.

(ii) By \cite[Proposition 3.3]{Wal}, any finite 2-dimensional
simplicial complex with free fundamental group has the homotopy type
of a wedge of 1-spheres and 2-spheres.
\qed

\section{The simplification of a complex}

Let $\H\; = (V,H)$ and $\H' = (V',H')$ be simplicial complexes. A {\em
  simplicial 
  map} from $\H$ to $\H'$ is a mapping $\p:V \to V'$ such that $X\p
\in H'$ for every $X \in H$ (that is, $\p$ sends simplices to
simplices). This simplicial map is {\em rank-preserving} if $|X\p| =
|X|$ for every $X \in H$.

Let $\H\; = (V,H) \in \BoR$. We define an equivalence relation $\eta_{\H}$
on $V$ by
$$a \eta_{\H} b \; \mbox{ if } \; \oo{a} = \oo{b}.$$
If no confusion arises, we omit the index from $\eta_{\H}$.

It follows from (\ref{cheta}) that $a \eta b$ if and only if $ab \notin H$. 
If $M$ is a boolean matrix representation of $\H$, it is easy to see
that $a \eta b$ if and only if
the column vectors $M[\,\underline{\hspace{.2cm}},a]$ and
$M[\,\underline{\hspace{.2cm}},b]$ are equal.
Indeed, $M[\,\underline{\hspace{.2cm}},a] = M[\,\underline{\hspace{.2cm}},b]$
implies $ab \notin H$ trivially and the converse follows from the fact
that there exist no zero columns in $M$ (since $P_1(V) \subseteq
H$). Note also that (\ref{cheta1}) implies that $\oo{p} = p\eta$ for
every $p \in V$. 

The following lemma enhances the role played by $\eta$ in the context
of rank-preserving simplicial maps.

\bl
\label{rseta}
Let $\H\; = (V,H) \in \BoR$ and let $\tau$ be an equivalence relation
on $V$. Then the following conditions are equivalent:
\bi
\item[(i)] $\tau$ is the kernel of some rank-preserving simplicial map
  $\p:\H \to \H'$ into some simplicial complex $\H'$;
\item[(ii)] $\tau \subseteq \eta_{\H}$.
\ei
\el

\proof
(i) $\Rw$ (ii). Let $a,b \in V$ and suppose that $(a,b) \notin
\eta$. Then $\oo{a} \neq \oo{b}$ and so $ab \in H$. Since $\p$ is a
rank-preserving simplicial map, it follows that $a\p \neq b\p$ and so
$(a,b) \notin \tau$. Thus $\tau \subseteq \eta$.

(ii) $\Rw$ (i). We define a
simplicial complex $\H/\tau = (V/\tau,H/\tau)$, where 
$$H/\tau = \{
\{ a_1\tau, \ldots,a_k\tau \} \mid a_1\ldots a_k \in H \}.$$
Let $\p:V \to V/\tau$ denote the canonical projection. By definition,
$\p$ is a simplicial map. We claim that
\beq
\label{retro1}
\mbox{
$\p$ is rank-preserving.}
\eeq
Indeed, every (nonempty) $X \in H$ admits an enumeration
$x_1,\ldots, x_k$ satisfying (\ref{derby3}) and so $\oo{x_i} \neq
\oo{x_j}$ whenever $i \neq j$. Thus 
$$x_i \tau x_j \Rw x_i\eta x_j \Rw \oo{x_i} = \oo{x_j}  \Rw i = j$$
and so $|X\p| = |X|$. Thus (\ref{retro1}) holds and so $\tau$ is
the kernel of some rank-preserving simplicial map.
\qed

Note that, if $\tau \subseteq \eta$, it follows from the
characterization of $H$ in (\ref{derby3}) that
\beq
\label{retro2}
\mbox{if $a_i\tau b_i$ for $i = 1,\ldots,k$, then 
$a_1\ldots a_k \in H$ if and only if $b_1\ldots b_k \in
H$.}
\eeq

We collect in the next result some of the properties of the simplicial
complexes $\H/\tau$ (using the notation introduced in the proof of
Lemma \ref{rseta}).

\bp
\label{retro}
Let $\H\; = (V,H) \in \BoR$ and let $\tau \subseteq \eta$ be an
equivalence relation on $V$. Let $\p:V \to
V/\tau$ denote the canonical projection. Then:
\bi
\item[(i)] ${\rm dim}(\H/\tau) = {\rm dim}\H$;
\item[(ii)] ${\rm Fl}\H\; = \{ F\p\inv \mid F \in {\rm Fl}(\H/\tau)
  \}$;
\item[(iii)] ${\rm Fl}\H\cong {\rm Fl}(\H/\tau)$;
\item[(iv)] $\H/\tau$ is boolean representable;
\item[(v)] $\H/\tau$ is simple if and only if $\tau = \eta$;
\item[(vi)] $\H$ is pure if and only if $\H/\tau$ is pure;
\item[(vii)] $\H$ is a matroid if and only if $\H/\tau$ is a matroid;
\item[(viii)] if $v, w \in V$ are such that $v\tau \neq w\tau$, then $v \edge w$
  is an edge of $\Gamma{\rm Fl}\H$ if and only if $v\tau \edge w\tau$ is
  an edge of $\Gamma{\rm Fl}(\H/\tau)$;
\item[(ix)] for every $X \subseteq V$, $$X \in {\rm fct}\H \mbox{ if
  and only if } (\p|_X \mbox{ is injective and }X\p \in {\rm fct}(\H/\tau)).$$ 
\item[(x)] if $\H/\tau$ is shellable, so is $\H$.
\ei
\ep

\proof
(i) It follows from the definition of $H/\tau$ and (\ref{retro1}).

(ii) Let $F \in {\rm Fl}(\H/\tau)$. Let $X\in H \cap 2^{F\p\inv}$ and
$p \in V \setminus F\p\inv$. Then $X\p \in (H/\tau) \cap 2^{F}$ and
$p\tau \in (V/\tau)\setminus F$, hence $F \in {\rm Fl}(\H/\tau)$ yields
$X\p \cup \{ p\tau \} \in H/\tau$. Since the elements of $X\p \cup \{
p\tau \}$ are all distinct, it follows easily from (\ref{retro2}) that
$X \cup \{ p \} \in H$. Thus $F\p\inv \in \flatx\H$.

To prove the opposite inclusion, we start by showing that
\beq
\label{retro4}
\mbox{if $Z \in \flatx\H$, then $Z\p \in
\flatx(\H/\tau)$.}
\eeq
Let $Y \in (H/\tau) \cap 2^{Z\p}$ and $p\tau \in (V/\tau) \setminus
(Z\p)$. We may write $Y = X\p$ for some $X \in H$.
Since $a\p\p\inv \subseteq \oo{a}$ for every $a \in V$, we
have $Z\p\p\inv \subseteq Z$. Hence $X \in H \cap 2^Z$. On the other
hand, $p\tau \in (V/\tau) \setminus (Z\p)$ implies $p \in
V \setminus Z$. Since $Z \in \flatx\H$, we get $X \cup \{ p \} \in H$
and so $Y \cup \{ p\tau \} \in H/\tau$. Therefore $Z\p \in
\flatx(\H/\tau)$ and so (\ref{retro4}) holds.

Let $Z \in \flatx\H$. Since we have already remarked that $Z\p\p\inv
\subseteq Z$ and the opposite inclusion holds trivially, we get
$Z = Z\p\p\inv \in \{ F\p\inv \mid F \in {\rm Fl}(\H/\tau) \}$.

(iii) By part (ii), the mapping
$$\begin{array}{rcl}
{\rm Fl}(\H/\tau)&\to&{\rm Fl}\H\\
F&\mapsto&F\p\inv
\end{array}$$
is bijective, and is clearly a poset isomorphism. Therefore it is a
lattice isomorphism.

(iv) Let $X \in H$ so that $X\p \in H/\tau$. In view of
(\ref{retro1}) and part (ii),
there exists some enumeration 
$x_1,\ldots,x_k$ of the elements of $X$ and some $F_0,\ldots,F_k \in
\flatx(H/\tau)$ such that 
$$F_0\p\inv \subset F_1\p\inv \subset \ldots \subset F_k\p\inv$$
and $x_i \in (F_i\p\inv)\setminus(F_{i-1}\p\inv)$ for $i =
1,\ldots,k$. It follows that $F_0 \subset \ldots \subset F_k$ and
$x_i\p \in F_i\setminus F_{i-1}$ for every $i$, hence $X\p$ is a 
transversal of the successive differences for a chain in
$\flatx(\H/\tau)$. Therefore $\H/\tau$ is boolean representable.

(v) Given $X \subseteq V$, let $\cl_{\tau}(X\p)$ denote the closure of
$X\p$ in $\H/\tau$. We show that
\beq
\label{retro3}
\cl_{\tau}(X\p) = \oo{X}\p.
\eeq

Indeed, by (\ref{retro4}) we have $\oo{X}\p \in \flatx(H/\tau)$, and
trivially $X\p \subseteq \oo{X}\p$. Suppose now that $F \in
\flatx(H/\tau)$ contains $X\p$. By part (ii), we have $X \subseteq
F\p\inv \in \flatx\H$, hence $\oo{X} \subseteq F\p\inv$ by minimality
and so $\oo{X}\p \subseteq F$. Therefore (\ref{retro3}) holds.

Suppose now that $(a,b) \in \eta \setminus \tau$. Then (\ref{retro3})
yields $\cl_{\tau}(a\p) = \oo{a}\p = \oo{b}\p = \cl_{\tau}(b\p)$ and
so $\{ a\tau,b\tau\} \notin H\tau$ by
(\ref{cheta}). Therefore $\H/\tau$ is not simple.  

Finally, assume that $\tau = \eta$. Let $a,b \in V$ be such that $a\eta
\neq b\eta$. Then $\oo{a} \neq \oo{b}$ and by (\ref{cheta}) we get $ab
\in H$. Hence $\{ a\eta,b\eta \} \in H/\eta$ and so $\H/\eta$ is simple. 

(vi) Considering transversals of successive differences, it is
immediate that a boolean representable simplicial complex is pure if
and only if its lattice of flats satisfies the Jordan-Dedekind
condition (all the maximal chains have the same length). Now we use
part (iii).

(vii) It is well known that $\H$ is a matroid if and only if
$\flatx\H$ is geometric \cite[Theorem 1.7.5]{Oxl}. Now we use
part (iii).

(viii) Assume that $v \edge w$ is an edge of $\Gamma\flatx\H$. By part
(ii), there exists some $F \in \flatx(\H/\tau)$ such that 
$vw \subseteq F\p\inv \subset V$. It follows that $\{ v\tau,w\tau\}
\subseteq F \subset V/\tau$, hence $v\tau \edge w\tau$ is an edge of
$\Gamma\flatx(\H/\tau)$. 

Conversely, assume that $v\tau \edge w\tau$ is an edge of
$\Gamma\flatx(\H/\tau)$. Then there exists some $F \in
\flatx(\H/\tau)$ such that  
$\{ v\tau,w\tau\} \subseteq F \subset V/\tau$. Hence $vw \subseteq F\p\inv
\subset V$. Since $F\p\inv \in \flatx\H$ by part (ii), it follows that 
$v \edge w$ is an edge of $\Gamma\flatx\H$.

(ix) Let $X \in \fct\H$. Then $X\p \in
H/\tau$ and $\p|_X$ is injective by (\ref{retro1}). Suppose that $X\p \subset Y$
for some $Y \in H/\tau$. We may write $Y = X\p \cup Z\p$ with
$Z$ minimal. It follows from the minimality of $Z$ that $\p|_{X\cup Z}$
is injective, hence 
$X \cup Z \in H$ in view of (\ref{retro2}), contradicting $X \in
\fct\H$. Therefore $X\p \in {\rm fct}(\H/\tau)$.

Conversely, assume that $\p|_X$ is injective and $X\p \in {\rm
  fct}(\H/\tau)$. In view of (\ref{retro2}), we have $X \in H$.
Suppose that $X \cup \{ p \} \in H$ with $p \in
V \setminus X$. By (\ref{retro1}), $\p|_{X\cup \{ p \}}$ is injective
and $(X \cup \{ p \})\p \in H/\tau$, hence $X\p \subset (X \cup \{ p \})\p
\in H$, contradicting $X\p \in \fct(\H/\tau)$. Therefore $X \in
\fct\H$ and the equivalence holds.

(x) We may assume that $|V| = |V/\tau| +1$, and then apply this case
successively. Assume that $\{ a_1,a_2\}$ is the only nonsingular
$\tau$-class.

Let $B_1,\ldots,B_t$ be a shelling of $\H/\tau$. For $k = 1,2$, let 
$\psi_k:V/\tau \to V$ be defined by
$$x\p\psi_k = \left\{
\begin{array}{ll}
a_k&\mbox{ if }x \in \{ a_1,a_2\}\\
x&\mbox{ otherwise}
\end{array}
\right.$$
Consider the sequence
\beq
\label{ers}
B_1\psi_1,B_1\psi_2,B_2\psi_1,B_2\psi_2,\ldots, B_t\psi_1,B_t\psi_2.
\eeq
We have $B_i\psi_1 = B_i\psi_2$ if and only if $a_1\p \notin B_i$. To avoid
repetitions, we remove from
(\ref{ers}) all the entries $B_i\psi_2$ such that $a\p \notin B_i$. We
refer to this sequence as {\em trimmed} (\ref{ers}).

It follows from part (ix) that trimmed
(\ref{ers}) is an enumeration of the facets of $\H$. We prove it is a shelling. 

Let $i \in \{ 2,\ldots,t\}$ and assume that $|B_i| \geq 2$. Write
$$I(B_{i}) = (\cup_{j=1}^{i-1} 2^{B_{j}}) \cap 2^{B_{i}}, \quad
I'(B_i\psi_1) = ((\cup_{j=1}^{i-1} 2^{B_{j}\psi_1})\cup (\cup_{j=1}^{i-1}
2^{B_{j}\psi_2})) \cap 2^{B_{i}\psi_1}.$$ 
It is immediate that $I'(B_i\psi_1) = (I(B_i))\psi_1$. Since
$B_1,\ldots,B_t$ is a shelling of 
$\H/\tau$, then $(B_i,I(B_i))$ is pure of dimension $|B_i|-2$. Thus 
$(B_i\psi_1,I'(B_i\psi_1))$ is pure of dimension $|B_i\psi_1|-2$. 

Assume now that $i \in \{ 1,\ldots,t\}$, $a_1\p \in B_i$ and $|B_i|
\geq 2$. Write 
$$I'(B_i\psi_2) = ((\cup_{j=1}^{i} 2^{B_{j}\psi_1})\cup (\cup_{j=1}^{i-1}
2^{B_{j}\psi_2})) \cap 2^{B_{i}\psi_2}.$$

Assume first that $i = 1$. Then 
$$I'(B_1\psi_2) = 2^{B_1 \setminus \{ a_2 \}},$$
hence 
$(B_1\psi_2,I'(B_1\psi_2))$ is pure of dimension $|B_1\psi_2|-2$.

Thus we may assume that $i > 1$. It is easy to check that
\beq
\label{retro5}
I'(B_i\psi_2) = (I(B_i) \cup 2^{B_i \setminus \{ a_1\p \}})\psi_2.
\eeq
Since $(B_i,I(B_i))$ is pure of dimension $|B_i|-2$, it follows that
$(B_i,I(B_i) \cup 2^{B_i \setminus \{ a_1\p \}})$ has also dimension
$|B_i|-2$. Since the only new facet with respect to $(B_i,I(B_i))$
is possibly $B_i \setminus \{ a_1\p\}$, then 
$(B_i,I(B_i) \cup 2^{B_i \setminus \{ a_1\p \}})$ is also pure. In
view of (\ref{retro5}),
$(B_i\psi_2,I'(B_i\psi_2))$
is pure of dimension $|B_i\psi_2|-2$. 
Therefore trimmed (\ref{ers}) is a shelling of $\H$
and we are done.
\qed

Part (ii) implies that the maps $\p$ constitute a particular case of
maps known in matroid theory as {\em strong maps} \cite[Chapter 8]{Whi2}.

We could not prove so far the converse of Proposition \ref{retro}(x),
which remains an open problem.
However, it follows from Theorem \ref{shelltwo} that it holds for the
particular case of $\eta$ and dimension 2.

From now on, and in view of part (v), we shall refer to $\H_S = \H/\eta$ as
the {\em simplification} of $\H$.

The next result shows how we can produce a boolean representation for
$\H_S$ from a boolean representation of $\H$.

\bp
\label{scon}
Let $M$ be an $R \times V$ boolean matrix representation of the
simplicial complex $\H\; = (V,H)$. Let $M'$ be the matrix obtained
from $M$ by removing repeated columns. Then $M'$ is a boolean matrix
representation of $\H_S$.
\ep

\proof
By the remark following the definition of $\eta$, we have $a\eta b$ if
and only if $M[\,\underline{\hspace{.2cm}},a] =
M[\,\underline{\hspace{.2cm}},b]$. Hence we may view the column space
of $M'$ as $V/\eta$. Let $\p:V \to V/\eta$ denote the canonical projection.

Let $X \in H$ so that $X\p \in H/\eta$. Then there exists some $Y
\subseteq R$ such that $M[Y,X]$ is nonsingular. Then $M[Y,X]$ has no
repeated columns and so $M'[Y,X\p]$ is nonsingular. Thus $X\p$ is
$M'$-independent. 

Conversely, assume that $X' \subseteq V/\eta$ is
$M'$-independent. Write $X' = X\p$ with $|X|$ minimum. Then there exists some $Y
\subseteq R$ such that $M'[Y,X']$ is nonsingular. Since $|X| = |X'|$
by minimality, it follows easily that $M[Y,X]$ and $M'[Y,X']$ have the
same structure, hence $M[Y,X]$ is nonsingular. Therefore $X \in H$ and
so $X' = X\p \in H/\eta$ as required.
\qed 

We end this section by discussing how the fundamental groups of $\H$
and $\H_S$ are related. 

\bp
\label{fhhs}
Let $\H$ be a boolean representable simplicial
complex of dimension $\geq 2$. Then the following conditions are
equivalent:
\bi
\item[(i)] $\pi_1(||\H||) \cong \pi_1(||\H_S||)$;
\item[(ii)] every $H$-trivial connected components of $\Gamma{\rm
  Fl}\H$ has size 1. 
\ei
\ep

\proof
We show that 
\beq
\label{fhhs1}
\mbox{$\Gamma{\rm Fl}\H$ and $\Gamma{\rm Fl}\H_S$ have the same number
  of connected components.}
\eeq

Let $\H\; = (V,H)$ and denote by $\p:V \to V/\eta$ the canonical
projection.
Let $C_1, \ldots, C_m \subseteq V$ denote the connected components
of $\Gamma\flatx\H$ and let $C'_1, \ldots, C'_n \subseteq V/\eta$ denote the
connected components of $\Gamma\flatx\H_S$.

Given $i \in \{ 1,\ldots,m\}$, it follows easily from 
Proposition \ref{retro}(viii) that 
$C_i\p \subseteq C'_{k_i}$ for some $k_i \in \{ 1,\ldots,n\}$. Since
$V/\eta = C_1\p \cup \ldots \cup C_m\p$, it follows that $m \geq n$.

Suppose now that $k_i = k_j$ for some distinct $i,j \in \{
1,\ldots,m \}$. Take vertices $v_i$ and $v_j$ in $C_i$ and $C_j$,
respectively. If $v_i\eta \neq v_j\eta$, it follows easily from 
Proposition \ref{retro}(vi)
that $v_i,v_j$ are connected by some path, a contradiction. Hence we
may assume that $v_i\eta v_j$ and so $\oo{v_i}
= \oo{v_j}$ in $\H$. 

But $\H_S$ is simple, hence $\{ v_i\eta\} \in \flatx\H_S$ and so
$v_i\p\p\inv \in \flatx\H$ by Proposition \ref{retro}(ii). Since
$\{ v_i\eta\}$ and $V/\eta$ are 
distinct flats of $\H_S$, it also follows from Proposition
\ref{retro}(ii) that $v_i\p\p\inv 
\neq (V/\eta)\p\inv = V$, hence $v_i \edge v_j$ should be an edge of
$\Gamma\flatx\H$, a contradiction. Thus the correspondence $i \mapsto
k_i$ is injective and so $m = n$.

Therefore $\Gamma\flatx\H$ and $\Gamma\flatx\H_S$ have the same number
of connected components. 


Assume that $\Gamma{\rm Fl}\H$ has $s$
$H$-nontrivial connected components and $r$ $H$-trivial connected
components of sizes $f_1,\ldots,f_r$. By Theorem \ref{fung},
$\pi_1(||\H||)$ is a free group of rank 
$$\binom{s-1}{2} +(s-1)(f_1+\ldots +f_r) + \sum_{1 \leq i < j \leq r} f_if_j.$$

On the other hand, in view of (\ref{fhhs1}) and Corollary
\ref{simfun}, $\pi_1(||\H_S||)$ is a free 
group of rank 
$$\binom{s+r-1}{2} = \frac{(s+r-1)(s+r-2)}{2} = \frac{(s-1)(s-2) +
  (2s-3)r + r^2}{2} = \binom{s-1}{2} + (s-1)r + \binom{r}{2}.$$
Now $(s-1)(f_1+\ldots +f_r) \geq$ and $\sum_{1 \leq i < j \leq r}
f_if_j \geq \binom{r}{2}$, and both equalities hold if and only if
$f_1 = \ldots = f_r = 1$.
\qed

The following is one of the simplest examples with $\pi_1(||\H||) \neq
\pi_1(||\H_S||)$. 

\be
\label{chhs}
Let $V = 12345$ and $H = (P_{\leq 2}(V) \setminus 45) \cup \{ 123,
124, 125\}$. Then $\H\; = (V,H)$ is a boolean representable simplicial
complex of dimension $\geq 2$ suvh that $\pi_1(||\H||) \not\cong
\pi_1(||\H_S||)$.
\ee

Indeed, it is easy to check that
$$\flatx\H\; = \{ \emptyset, 1,2,3, 12, 45, V\}$$
and $\H$ is a boolean representable. Its graph of flats is
$$1 \edge 2 \hspace{1cm} 3 \hspace{1cm} 4 \edge 5,$$
hence the $H$-trivial connected components of $\Gamma{\rm Fl}\H$ have
size 1 and 2, respectively. Now the claim follows from Proposition
\ref{fhhs}.

\medskip

Note that there is a natural embedding of $\pi_1(||\H_S||)$ into
$\pi_1(||\H||)$ (since $\H_S$ is isomorphic to a restriction of $\H$ to a
cross-section of $\eta$) and this embedding splits since
$\pi_1(||\H_S||)$ is a free factor of $\pi_1(||\H||)$.

\section{Shellability and sequentially Cohen-Macaulay in dimension 2}

We discuss in this section shellability for boolean representable
simplicial complexes of dimension 2. The simple case was completely
solved in \cite[Theorem 7.2.8]{RSm}, now we generalize this theorem to arbitrary
boolean representable simplicial complexes of dimension 2. 

We consider also another property of topological significance, 
sequentially Cohen-Macaulay. It is often associated with shellability
since a shellable complex is necessarily sequentially Cohen-Macaulay
\cite{BWW,Sta2}. We need to introduce a few concepts and
notation before defining it.


Assume that $\dim\H\; = d$.
For $m = 0,\ldots,d$, we define the complex 
$\pure_m(\H) = (V_m,H_m)$ to be the subcomplex of $\H$ generated by
all the faces of $\H$ of dimension $m$. Clearly, $\pure_m(\H)$ is the
largest pure subcomplex of $\H$ of dimension $m$. 

In view of \cite[Theorem 3.3]{Duv}, we say that $\H$ is {\em
  sequentially Cohen-Macaulay} if
$$\tilde{H}_k(\pure_m(\lk(X))) = 0$$
for all $X \in H$ and $k < m \leq d$.



We start with the following lemma.

\bl
\label{scm}
Let $\H$ be a sequentially Cohen-Macaulay simplicial complex of
dimension 2. Then the simplification $\H_S$ is sequentially Cohen-Macaulay.
\el

\proof
Write $\H\; = (V,H)$.
Since $\dim\H_S = 2$ by Proposition \ref{retro}(i), we have to prove
the following facts:
\bi
\item[(1)] $\pure_2(\H_S)$ is connected;
\item[(2)] $\pure_1(\H_S)$ is connected;
\item[(3)] $\pure_1(\lk(v\eta))$ is connected for every $v \in V$;
\item[(4)] $\tilde{H}_1(\pure_2(\H_S)) = 0$.
\ei
We assume of course the similar statements for $\H$.

(1) Let $a\eta,b\eta$ denote two distinct vertices from
$\pure_2(\H_S)$. Then there exist $\{ a\eta,a'\eta,a''\eta\}$, $\{
b\eta,b'\eta,b''\eta\} \in (H/\eta) \cap P_3(V/\eta)$. In view of
(\ref{retro2}), we have $aa'a'',bb'b'' \in H \cap P_3(V)$, hence $a,b$
are two distinct vertices from $\pure_2(\H)$. Since $\pure_2(\H)$ is
connected, there exists in $\pure_2(\H)$ a path of the form
$$a = c_0 \edge c_1 \edge \ldots \edge c_n = b$$
for some $n \geq 1$. Let $i \in \{ 1,\ldots,n\}$. Since $c_{i-1}c_i$
is an edge of $\pure_2(\H)$,
there exists some $c'_i$ such that $c_{i-1}c_ic'_i \in H \cap
P_3(V)$. In view of (\ref{retro1}), we get $\{
c_{i-1}\eta,c_i\eta,c'_i\eta\} \in (H/\eta) \cap P_3(V/\eta)$. It
follows that
$$a\eta = c_0\eta \edge c_1\eta \edge \ldots \edge c_n\eta = b\eta$$
is a path in $\pure_2(\H_S)$ and so $\pure_2(\H_S)$ is connected.

(2) Similar to (1).

(3) Let $a\eta,b\eta$ denote two distinct vertices from
$\pure_1(\lk(v\eta))$. Then there exist
some edges $a\eta \edge a'\eta$, $b\eta \edge b'\eta$ in
$\lk(v\eta)$. Hence $\{ a\eta,a'\eta,v\eta\}, \{
b\eta,b'\eta,v\eta\} \in (H/\eta) \cap P_3(V/\eta)$. By
(\ref{retro2}), we get $aa'v,bb'v \in H \cap P_3(V)$, hence 
$a \edge a'$ and $b \edge b'$ are edges in $\lk(v)$ and so $a,b$
are two distinct vertices from $\pure_1(\lk(v))$. Since $\pure_1(\lk(v))$ is
connected, there exists in $\pure_1(\lk(v))$ a path of the form
$$a = c_0 \edge c_1 \edge \ldots \edge c_n = b$$
for some $n \geq 1$. Let $i \in \{ 1,\ldots,n\}$. Since $c_{i-1}c_i$
is an edge of $\pure_1(\lk(v))$, we have $c_{i-1}c_iv \in H \cap
P_3(V)$ and so (\ref{retro1}) yields
$\{ c_{i-1}\eta,c_i\eta,v\eta\} \in (H/\eta) \cap P_3(V/\eta)$. It
follows that
$$a\eta = c_0\eta \edge c_1\eta \edge \ldots \edge c_n\eta = b\eta$$
is a path in $\pure_1(\lk(v\eta))$ and so $\pure_1(\lk(v\eta))$ is connected.

(4) Fix a cross section $V_0 \subseteq V$
for $\eta$. We consider the ordering of $V/\eta$ induced by the restriction
of the ordering of $V$ to $V_0$.

Suppose that $\tilde{H}_1(\pure_2(\H_S)) \neq 0$. Let 
$\partial_k$ (respectively $\partial'_k$) denote the $k$th boundary
map of $\pure_2(\H)$ (respectively $\pure_2(\H_S)$). Since
$\ker\partial'_1 / \im\partial'_{2} = \tilde{H}_1(\pure_2(\H_S)) \neq
0$, there exist some distinct edges $X_1,\ldots, X_m$ in
$\pure_2(\H_S)$ and some $n_1,\ldots,n_m \in \Z$ such that
$\sum_{i=1}^m n_iX_i \in \ker\partial'_1 \setminus
\im\partial'_{2}$. 
Write $X_i = \{ a_i\eta,b_i\eta\}$ with $a_i,b_i \in V_0$ and $a_i <
b_i$. By definition 
of $\pure_2(\H_S)$, there exists some $c_i \in V_0$ such that $\{
a_i\eta,b_i\eta,c_i\eta\} \in (H/\eta) \cap P_3(V/\eta)$. In view of
(\ref{retro2}), we have $a_ib_ic_i \in H \cap P_3(V_0)$, hence $a_ib_i$
is an edge from $\pure_2(\H)$. Now
$$0 = (\sum_{i=1}^m n_iX_i)\partial'_1 = \sum_{i=1}^m n_i(b_i\eta-a_i\eta)$$
yields $\sum_{i=1}^m n_i(b_i-a_i) = 0$ since $V_0$ is a cross-section
for $\eta$ and so $\sum_{i=1}^m n_i(a_ib_i) \in \ker\partial_1$. 

Since $0 = \tilde{H}_1(\pure_2(\H)) = \ker\partial_1 /
\im\partial_{2}$, we must have 
\beq
\label{scm1}
\sum_{i=1}^m n_i(a_ib_i) = (\sum_{j=1}^r k_j(x_jy_jz_j))\partial_2
\eeq
for some distinct triangles $x_jy_jz_j$ in
$\pure_2(\H)$ and $k_j \in \Z$. 
Since $a_i,b_i \in V_0$ for every $i$,
we may assume 
that $x_j < y_j < z_j$ and $x_j,y_j,z_j \in V_0$ for every $j$:
indeed, we may replace 
each letter in $V\setminus V_0$ by its representative in $V_0$, and
remain inside $\pure_2(\H)$ by (\ref{retro2}). In view of
(\ref{retro1}), $\{ x_j\eta, y_j\eta, 
z_j\eta \}$ is a triangle in $\H_S$ (and therefore in $\pure_2(\H_S)$)
for $j = 1,\ldots,r$. Now (\ref{scm1}) yields
$$\sum_{i=1}^m n_i(a_ib_i) = \sum_{j=1}^r k_j(y_jz_j - x_jz_j +
x_jy_j)$$ 
and consequently
$$\sum_{i=1}^m n_i\{ a_i\eta, b_i\eta \} = \sum_{j=1}^r k_j(\{
y_j\eta, z_j\eta\} - \{ x_j\eta, z_j\eta \} + \{ x_j\eta, y_j\eta \}
).$$
Since $x_j\eta < y_j\eta < z_j\eta$, we get 
$$\sum_{i=1}^m n_iX_i = (\sum_{j=1}^r k_j\{ x_j\eta, y_j\eta,
z_j\eta\})\partial'_2 \in \im\partial'_{2},$$
a contradiction. Therefore $\tilde{H}_1(\pure_2(\H_S)) = 0$ as
required.
\qed
 
We may now prove one of our main theorems. The simple case (for
dimension 2) had been established in \cite[Corollary 7.2.9]{RSm}.

\bt
\label{shelltwo}
Let $\H$ be a boolean representable simplicial complex of dimension 2.
Then the following conditions are equivalent:
\bi
\item[(i)] $\H$ is shellable;
\item[(ii)] $\H$ is sequentially Cohen-Macaulay;
\item[(iii)] $\Gamma${\rm Fl}$\H_S$ contains at most two connected
components or contains exactly one nontrivial connected
component.
\ei
\et

\proof
(i) $\Rw$ (ii). By \cite{BWW,Sta2}.

(ii) $\Rw$ (i). By Lemma \ref{scm}, $\H_S$ is sequentially
Cohen-Macaulay. It follows from \cite[Corollary 7.2.9]{RSm} that
$\H_S$ shellable. Therefore $\H$ is shellable by Proposition \ref{retro}(x).

(i) $\Rw$ (iii). 
We adapt the proof of \cite[Lemma 7.2.7]{RSm}.

Let $C_1,\ldots,C_m$ denote the connected components of $\Gamma${\rm
  Fl}$\H_S$. We suppose that $m \geq 3$ and at least $C_1,C_2$ are
nontrivial. Since $\H_S$ is simple of dimension 2, we know by
\cite[Lemma 6.4.3]{RSm} that
\beq
\label{shelltwo1}
\mbox{if $pqr \in P_3(V/\eta)$, $p \edge q$ is an edge of
  $\Gamma\flatx\H_S$ but $p \edge r$ is not, then }pqr \in H/\eta.
\eeq
 
Let $\p:V \to V/\eta$ be the canonical projection.
For $i = 1,\ldots,m$, let
$V_i = \{ v \in V \mid v\p \in C_i\}$. It follows that $V = V_1 \cup
\ldots \cup V_m$ constitutes a partition of $V$. We show that
\beq
\label{shelltwo2}
\mbox{if $pqr \in H \cap P_3(V)$, then }p,q,r
\mbox{ belong to at most two distinct }V_i.
\eeq

It follows from (\ref{retro1}) that
$\{ p\p, q\p, r\p\} \in (H/\eta) \cap P_3(V/\tau)$. By Proposition
\ref{retro}, $\H_S$ is a simple boolean
representable simplicial complex of dimension 2, so it follows from
\cite[Lemma 6.4.4]{RSm} that the three 
vertices $p\p,q\p,r\p$ belong to at most two connected
components of $\Gamma\flatx\H_S$. Therefore (\ref{shelltwo2}) holds.

We split now the discussion into two cases.
Suppose first that $\Gamma${\rm Fl}$\H_S$ has a trivial connected
component $C_k$. Let $v$ be its single vertex. We consider the link $\lk(v)$.
By \cite{BW} (see also \cite[Proposition 7.1.5]{RSm}), $\H$ shellable
implies $\lk(v)$ shellable. Let $p_i\eta \edge q_i\eta$ be an edge of
$C_i$ for $i = 1,2$. By (\ref{shelltwo1}), we have $\{ p_i\eta,
q_i\eta, v\eta \} \in H/\tau$. By (\ref{retro2}), we get $p_iq_iv \in
H$, hence $p_iq_i \in H/v$ and so $\lk(v)$ has dimension 1.

The facets of a complex of dimension 1 are the edges and the isolated
vertices. It is immediate that such a 
complex is shellable if and only the complex has a unique nontrivial
connected component. Therefore, since $\lk(v)$ is shellable of dimension
1, the edges $p_1q_1,p_2q_2 \in H/v$ must belong to the same
connected component of $\lk(v)$. Hence there exist distinct $r_0,\ldots,r_n \in
V\setminus \{ v \}$ such that $r_0 \in p_1q_1$, $r_n \in p_2q_2$ and
$r_{j-1}r_j \in H/v$ for $j = 1,\ldots,n$. 

Now we have $r_{j-1}r_jv \in H$. Since $v$ is an isolated vertex of
$\Gamma\flatx\H_S$, then $H \cap P_2(V_k) = \emptyset$ by
(\ref{retro1}). Hence (\ref{shelltwo2}) yields 
$r_{j-1},r_j \in V_i$ for some $i \in \{ 1,\ldots,m \} \setminus \{ k
\}$. Thus $r_0, r_n \in V_i$. But $r_0 \in p_1q_1$ and $r_n \in
p_2q_2$ imply $r_0\in C_1$ and $r_n \in C_2$, a contradiction.  

Therefore we may assume that all the connected components $C_1,\ldots,
C_m$ of $\Gamma\flatx\H_S$ are nontrivial. 

Suppose that $pq \in
H \cap P_2(V)$. By (\ref{cheta}), we have $p\eta \neq
q\eta$. If $p\eta \edge q\eta$ is an edge of $\Gamma\flatx\H_S$, let
$r \in V$ be such that $r\eta \notin \oo{\{ p\eta,q\eta\}}$. Then $\{
p\eta,q\eta, r\eta\} \in H/\eta$ and in view of (\ref{retro2}) we get
$pqr \in H \cap P_3(V)$. Thus $\H$ has
no 1-facets.

On the other hand, given $p \in V$, we may take $q \in V \setminus
p\p\inv$. Since $\H_S$ is simple, we have $\{ p\eta,q\eta \} \in H/\eta$, 
yielding $pq \in H$ in view of (\ref{retro2}). Therefore
every facet of $\H$ has dimension $2$.

Let $B_1,\ldots,B_t$ be a shelling of $\H$. For $k = 1,\ldots,t$,
define a graph $\Gamma_k = (W_k,E_k)$ by
$$W_k = \cup_{j=1}^k B_j,\quad E_k = \cup_{j=1}^k P_2(B_j).$$
It follows easily from the definition of shelling that each $\Gamma_k$
is connected.

We say that $p,q \in W_k$ have the same color if $p,q \in V_i$ for
some $i \in \{ 1,\ldots,m\}$. We write $p\gamma_k q$ if there exists a
monochromatic path of the form
$$p = r_0 \edge r_1 \edge \ldots \edge r_n = q$$
in $\Gamma_k$ for some $n \geq 0$. It is immediate that $\gamma_k$ is
an equivalence 
relation on $W_k$. We define a graph $\oo{\Gamma_k} =
(\oo{W_k},\oo{E_k})$ by taking 
$\oo{W_k} = \{ p\gamma_k \mid p \in W_k\}$ and 
$$\oo{E_k} = \{ \{ p\gamma_k, q\gamma_k\} \mid p\gamma_k \neq
q\gamma_k \mbox{ and } pq \in E_k\}.$$
We prove that 
\beq
\label{shelltwo3}
\oo{\Gamma_k}\mbox{ is a tree for }k = 1,\ldots,t
\eeq
by induction on $k$.

In view of (\ref{shelltwo2}), $\oo{\Gamma_1}$ has at most two
vertices, hence a tree. Assume now that $k > 1$ and
$\oo{\Gamma_{k-1}}$ is a tree. We consider several cases and subcases:
\smallskip

\noindent
\underline{Case 1}: $B_k \not\subseteq W_{k-1}$.
\smallskip

\noindent
Since $B_k$ has dimension 2, then $(B_k,I(B_k))$ is pure of dimension
1, hence we may write $B_k = pqr$ with $pq \in E_{k-1}$ and $r \notin
W_{k-1}$. By (\ref{shelltwo2}), the vertices $p,q,r$ have at most two
different colors.
\smallskip

\noindent
\underline{Subcase 1.1}: $r$ has the same color as $p$ or $q$.
\smallskip

\noindent
Then $\oo{\Gamma_k} = \oo{\Gamma_{k-1}}$, hence a tree by the
induction hypothesis.
\smallskip

\noindent
\underline{Subcase 1.2}: $r$ has a different color from $p$ and $q$.
\smallskip

\noindent
Then $p \gamma_{k-1} q$ and so $\oo{\Gamma_k}$ is obtained from
$\oo{\Gamma_{k-1}}$ by adjoining the edge $p\gamma_{k-1} = p\gamma_k \edge
r\gamma_k$. Since $\oo{\Gamma_{k-1}}$ is a tree, $\oo{\Gamma_{k}}$ is
a tree as well.

\smallskip

\noindent
\underline{Case 2}: $B_k \subseteq W_{k-1}$.
\smallskip

\noindent
We may assume that $E_{k-1} \subset E_k$.
Since $B_k$ has dimension 2, then $(B_k,I(B_k))$ is pure of dimension
1, hence we may write $B_k = pqr$ with $pq,qr \in E_{k-1}$ and $pr \notin
E_{k-1}$. By (\ref{shelltwo2}), the vertices $p,q,r$ have at most two
different colors.
\smallskip

\noindent
\underline{Subcase 2.1}: $q$ has the same color as $p$ or $r$.
\smallskip

\noindent
Then $\oo{\Gamma_k} = \oo{\Gamma_{k-1}}$, hence a tree by the
induction hypothesis.
\smallskip

\noindent
\underline{Subcase 2.2}: $q$ has a different color from $p$ and $r$.
\smallskip

\noindent
Then $p$ and $r$ have the same color. If $p \gamma_{k-1} r$, then 
$\oo{\Gamma_k} = \oo{\Gamma_{k-1}}$, hence we may assume that $(p,r)
\notin \gamma_{k-1}$. It follows that $\oo{\Gamma_k}$ is obtained from
$\oo{\Gamma_{k-1}}$ by identifying the (non adjacent) vertices $p\gamma_{k-1}$ and
$q\gamma_{k-1}$. It is well known that folding such a pair of adjacent
edges in a tree still yields a tree.



Therefore $\oo{\Gamma_k}$ is a tree in all cases and so
(\ref{shelltwo3}) holds.

Let $p_i \edge q_i$ be an edge in $C_i$ for $i = 1,2$ and let $v$ be a
vertex in $C_3$. By (\ref{shelltwo1}), we have $p_1q_1p_2, p_1q_1v,
p_2q_2v \in H$. Since all the facets in $\H$ have dimension 2, we have
$E_t = H \cap P_2(V)$, hence
$$
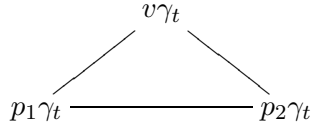
\xymatrix{
&{v\gamma_t} \ar@{-}[dl] \ar@{-}[dr] &\\
{p_1\gamma_t} \ar@{-}[rr] && {p_2\gamma_t}
}$$
is a triangle in $\oo{\Gamma_t}$, contradicting (\ref{shelltwo3}).
Therefore condition (ii) must hold.

(iii) $\Rw$ (i). By \cite[Theorem 7.2.8]{RSm}, $\H_S$ is shellable,
  which implies $\H$ shellable by Proposition \ref{retro}(x).
\qed

It is well known that a shellable simplicial complex has the
homotopy type of a wedge of spheres \cite{BW}. But in the case of BRSCs of
dimension 2, we already know from Theorem \ref{wed}(ii) that this is
always the case, despite there being such complexes that are not
shellable (see e.g. Example \ref{occur} for $t \geq 3$).

\section{The order complex of a lattice and EL-labelings}

Given a lattice $L$, let $C_{L^*}$ denote the set of totally ordered
subsets of $L^* = L \setminus \{ 0,1\}$ (chains). The {\em order
  complex} of $L$ is the 
simplicial complex $\ord(L) = (L^*,C_{L^*})$.

The concept of EL-labeling provides a famous sufficient condition for
shellability of the order complex of a lattice. 
Let $L$ be a lattice
and let $EH_L$ denote the set of edges in the Hasse diagram of
$L$. More formally, we can define $EH_L$ as the set of all ordered
pairs $(a,b) \in L \times L$ such that $b$ covers $a$ in $L$. 
Let $P$ be
a poset and let $\xi:EH_L \to P$ be a mapping. Given a maximal chain
$\gamma:\ell_0 < \ell_1 < \ldots < \ell_n$ in $L$ (so that
$(\ell_{i-1},\ell_{i}) 
\in EH_L$ for $i = 1,\ldots,n$), we define a word
$\gamma\xi$ on the alphabet $P$ by $\gamma\xi =
(\ell_0,\ell_{1})\xi\ldots(\ell_{n-1},\ell_{n})\xi$. The chain $\gamma$
is {\em increasing} if $(\ell_0,\ell_{1})\xi < \ldots
<(\ell_{n-1},\ell_{n})\xi$. Given $a,b \in L$ with $a < b$, we denote by
$[a,b]$ the subsemilattice of $L$ consisting of all $c \in L$
satisfying $a \leq c \leq b$. Clearly, $\xi:EH_L \to P$ induces also a
mapping on the maximal chains of $[a,b]$. Consider the lexicographic
ordering on $P^+$. We say that $\xi:EH_L \to P$ is an {\em
  EL-labeling} of $L$ if, for all $a,b \in L$ such that $a < b$:
\bi
\item
there exists a unique maximal chain $\gamma_0$ in $[a,b]$ such that
$\gamma\xi$ is increasing;
\item
$\gamma_0\xi < \gamma\xi$ for every other maximal chain $\gamma$ in $[a,b]$.
\ei

A fundamental theorem of Bj\"orner \cite{Bjo} states that if a lattice
$L$ admits an EL-labeling, then $\ord(L)$ is shellable. Moreover, it
is known that every
semimodular lattice admits an EL-labeling \cite[Exercise
  3.2.14(d)]{Wac}. In the case of boolean representable simplicial
complexes, the lattice of flats is semimodular if and only if the
complex is a matroid \cite[Theorem 1.7.5]{Oxl}.

The next result shows how a shelling of the order complex 
can provide a shelling of the original complex itself.

\bt
\label{ella}
Let $\H$ be a boolean representable simplicial complex. If
the order complex of ${\rm Fl}\H$ is shellable, so is $\H$.
\et

\proof
Write $L = \flatx\H$ and let $d = \dim\H = \dim(\ord(L)) +1$. The facets
of $\ord(L)$ can be identified (recall that we are looking at chains
in $L^*$ in $\ord(L)$) with the maximal chains in $L$,
i.e. subsets of $L$ of the 
form $B = \{ F_0,\ldots,F_{n}\}$ with 
\beq
\label{ella1}
\emptyset = F_0 \subset F_1 \subset \ldots \subset F_{n} = V
\eeq
and no intermediate flat $F_{i-1} \subset F' \subset F_i$ for $i =
1,\ldots, n$. Note that $n \leq d+1$. 
We define $B\tau$ to be the set of transversals of the maximal chain
(\ref{ella1}), i.e. $B\tau$ consists of all the subsets $\{ a_1,\ldots,
a_{n}\} \in P_{n}(V)$ such that $a_i \in F_i\setminus 
F_{i-1}$ for $i = 1,\ldots,n$. Note that $F_i = \oo{F_{i-1} \cup \{
  a_i\}}$ by maximality of (\ref{ella1}).

Assume that $B_1,\ldots,B_t$ is a shelling of $\ord(L)$. Then
$$\fct\H = \bigcup_{i=1}^t B_i\tau.$$
We intend to concatenate successive enumerations of
$B_1\tau,\ldots,B_t\tau$ so that, after removing repetitions, we get a
shelling of $\H$. 

We start with $B_1\tau$. Assuming that $B_1\tau$ is the set of
transversals of the chain (\ref{ella1}), we fix a total ordering
$<_1$ of
$V$ such that $a <_1 b$ whenever $a \in F_i\setminus F_{i-1}$, $b \in
F_j\setminus F_{j-1}$ and $i < j$. We may associate to each $B'_k
\in B_1\tau$ a (unique) word $a_1\ldots a_{n} \in V^{n}$ 
such that $B'_k = \{ a_1,\ldots, a_{n} \}$ and $a_i \in F_i\setminus
F_{i-1}$ for $i = 1,\ldots,n$. Then we order the elements of
$B_1\tau$ according to the lexicographical ordering of the associated
words.

Let us check the shelling condition for the facets in $B_1\tau$,
enumerated as $B'_1,\ldots,B'_p$. Let $k \in \{ 2,\ldots, p\}$. Let
$A \in I(B'_k)$. Then $B'_k$ is not the minimum facet (for the
lexicographic order) containing $A$. Hence there exists some $i \in \{
1,\ldots,n\}$ and some letters $b,c \in F_i \setminus F_{i-1}$ such
that $b <_1 c \in B'_k \setminus A$. It follows that $(B'_k \setminus \{
c \}) \cup \{ b \} = B'_j$ for some $j < k$ and so $A \subseteq B'_k
\setminus \{ c \} \in I(B'_k)$. Thus $(B'_k,I(B'_k))$ is pure of dimension
$n-2$.

Assume now that $j \in \{ 2,\ldots,t\}$ and we have already defined
enumerations for the facets in 
$B_1\tau \cup \ldots \cup B_{j-1}\tau$ so that the shelling condition is
satisfied. We may assume that $B_j\tau$ is the set of
transversals of the chain (\ref{ella1}). We fix a total ordering $<_j$ of
$V$ such that $a <_j b$ whenever $a \in F_i\setminus F_{i-1}$, $b \in
F_r\setminus F_{r-1}$ and $i < r$. Similarly to the case $j = 1$, we
associate to each $B'_k 
\in B_j\tau$ a (unique) word $a_1\ldots a_{n} \in V^{n}$ 
such that $B'_k = \{ a_1,\ldots, a_{n} \}$ and $a_i \in F_i\setminus
F_{i-1}$ for $i = 1,\ldots,n$. Then we order the elements of
$B_j\tau$ according to the lexicographical ordering of the associated
words, and we concatenate the new elements, say $B'_1,\ldots,B'_p$, to
the enumeration of the 
elements of $B_1\tau \cup \ldots \cup B_{j-1}\tau$ previously defined.

Assume that $q \in \{ 1,\ldots,p\}$ and $B'_q = \{ a_1,\ldots, a_{n}
\}$, where $a_i \in 
F_i\setminus F_{i-1}$ for $i = 1,\ldots,n$. Let $A \in I(B'_q)$, say
$A = \{ a_{u_1},\ldots, a_{u_s} \}$. Let $\wt{A} = \{ F_{u_1},\ldots,
F_{u_s} \} \in I(B_j)$. Since $(B_j,I(B_j))$ is
pure of dimension $n-2$, there exists some $j' < j$ such that $\wt{A}
\subseteq B_{j'}$ and $B_{j'}$ contains all the elements of $B_j$ but
one, say $F_i$. We may then assume that 
$B_{j'}$ originates from the chain
\beq
\label{ella2}
\emptyset = F_0 \subset \ldots \subset F_{i-1} \subset G_1 \subset
\ldots \subset G_w \subset
F_{i+1} \subset \ldots \subset F_{n} = V
\eeq
in $L$. Note that the $G_i$ must appear consecutively as a replacement
of the missing $F_i$ by maximality of
(\ref{ella1}). We claim that $B'_q \setminus \{ a_i\}$ is a partial
transversal of (\ref{ella2}) containing $A$.

Suppose that $a_i \in A$. Then $F_i \in \wt{A} \subseteq B_{j'}$, a
contradiction since (\ref{ella1}) is maximal and different from
(\ref{ella2}). Hence $a_i \notin A$ and so $A \subseteq  
B'_q \setminus \{ a_i\}$. To show that $B'_q \setminus \{ a_i\}$ is a partial
transversal of (\ref{ella2}), it is enough to note that 
$$a_{i+1} \in
F_{i+1}\setminus F_i \subseteq (F_{i+1}\setminus G_w)\cup \ldots \cup
(G_2 \setminus G_1) \cup (G_1 \setminus
F_{i-1}).$$
Thus $A \subseteq B'_q \setminus \{ a_i\} \in I(B'_q)$ and
so $(B'_q,I(B'_q))$ is pure of dimension $n-2$. By double induction on 
$q$ and $j$, this validates our construction of a shelling of $\H$.
\qed

The next example shows that the converse of Theorem \ref{ella} does
not hold.

\be
\label{noel}
Let $V = \{ 1,\ldots,6\}$ and let $\Gamma$ be the graph
$$\xymatrix{
1 \ar@{-}[r] & 2 \ar@{-}[r] & 3 \ar@{-}[r] & 4 & 5 \ar@{-}[r] & 6
}$$
Let 
$$H = P_{\leq 2}(V) \cup \{ X \in P_3(V) \mid \mbox{ at least two
  vertices in $X$ are adjacent in }\Gamma \}$$
and $\H\; = (V,H)$.
Then $\H$ is a shellable pure boolean representable simplicial complex
but the order complex of ${\rm Fl}\H$ is not shellable. 
\ee

Since there exist no isolated vertices in $\Gamma$, $\H$ is pure.
It is easy to compute the flats of $\H$, we have 
$$\flatx\H \; = P_{\leq 1}(V) \cup \{ 12, 23, 34, 56, V\}.$$
It follows easily that $\H$ is boolean representable. Moreover,
$\Gamma$ is indeed the graph of flats of $\H$, hence $\H$ is shellable
by Theorem \ref{shelltwo}. A possible shelling is
$$123, 124, 125, 126, 134, 156, 234, 235, 236, 256, 345, 346, 356, 456.$$

Now the facets of $\ord(\flatx\H)$ are 
$$\{ 1, 12\}, \{ 2,12\}, \{ 3,34\}, \{ 4,34\}, \{ 5,56\}, \{ 6,56\}.$$
It is well known that a graph is shellable if and only if has at most
one nontrivial connected component, hence $\ord(\flatx\H)$ is not shellable.



\bigskip

In the matroid case, we can combine Theorem \ref{ella} with the
aforementioned results of Bj\"orner on EL-labelings
to produce shellings for matroids (see \cite{Bjo}). Example \ref{noel} 
provides an example of a shellable pure boolean representable
simplicial complex which admits no EL-labeling of the lattice of
flats (otherwise $\ord(\flatx\H)$ would be shellable). Of course, this
simplicial complex is not a matroid.
The next example shows that the existence of EL-labelings is not
exclusive of matroids.

\be
\label{yesel}
Let $V = \{ 1,\ldots,7\}$ and let $\Gamma$ be the graph
$$\xymatrix{
1 \ar@{-}[r] & 2 \ar@{-}[r] & 3 \ar@{-}[r] & 4 \ar@{-}[r] & 5 \ar@{-}[r] & 
6 \ar@{-}[r] & 7
}$$
Let 
$$H = P_{\leq 2}(V) \cup \{ X \in P_3(V) \mid \mbox{ at least two
  vertices in $X$ are adjacent in }\Gamma \}$$
and $\H\; = (V,H)$.
Then $\H$ is a shellable pure boolean representable simplicial complex
which is not a matroid
and ${\rm Fl}\H$ admits an EL-labeling.
\ee

Since there exist no isolated vertices in $\Gamma$, $\H$ is pure.
It is easy to compute the flats of $\H$, we have 
$$\flatx\H \; = P_{\leq 1}(V) \cup \{ 12, 23, 34, 45, 56, 67, V\}.$$
It is easy to check now that $\H$ is boolean representable and
$\Gamma$ is indeed the graph of flats of $\H$. Thus $\H$ is shellable 
by Theorem \ref{shelltwo}.

The exchange property fails for $123$ and $57$, hence $\H$ is not a
matroid. The following diagram describes an EL-labeling
$\xi:EH_{\flatx\H} \to \mathbb{N}$. where the naturals are endowed
with the usual ordering.
$$
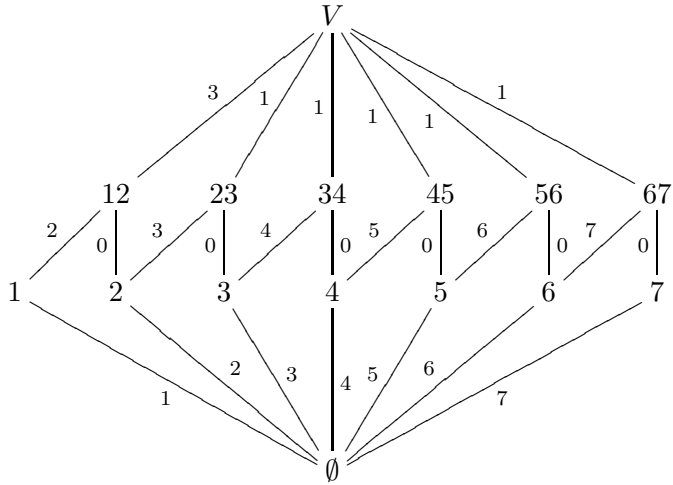
\xymatrix{
&&&V \ar@{-}[ddll]_3 \ar@{-}[dd]_1 \ar@{-}[ddl]_1 \ar@{-}[ddr]_1 \ar@{-}[ddrr]_1
  \ar@{-}[ddrrr]^1 &&&\\ &&&&&&\\
&12 &23 &34 &45 & 56 & 67\\
1 \ar@{-}[ur]^2 & 2 \ar@{-}[u]^0 \ar@{-}[ur]^3 & 3 \ar@{-}[u]^0
\ar@{-}[ur]^4 & 4 \ar@{-}[u]_0 \ar@{-}[ur]^5 & 5 \ar@{-}[u]^0
\ar@{-}[ur]^6 & 6 \ar@{-}[u]_0 \ar@{-}[ur]^7 & 7 \ar@{-}[u]^0\\
&&&&&&\\
&&&\emptyset \ar@{-}[uulll]^1 \ar@{-}[uull]_2 \ar@{-}[uul]_3 \ar@{-}[uu]_4
\ar@{-}[uur]^5 
\ar@{-}[uurr]^6 \ar@{-}[uurrr]_7 &&&
}$$

\section{Computing the flats}

In this section, we discuss the computation of the flats for a boolean
representable simplicial complex of fixed dimension $d$, and relate
these computations to the main results of the paper. The case $d \leq
1$ is straightforward and shall be omitted in most results.

We recall the $O$ notation from complexity theory. Let $P$ be an
algorithm defined for instances depending on parameters
$n_1,\ldots,n_k$. If $\p:\N^k \to \N$ is a function, we write $P \in
O((n_1,\ldots,n_k)\p)$ if there exist constants $K,L > 0$ such that
$P$ processes each instance of type $(n_1,\ldots,n_k)$ in time $\leq
K((n_1,\ldots,n_k)\p) + L$ (where time is measured as the number of
elementary operations performed).

Clearly, boolean matrices provide the most natural means of defining
a boolean representable simplicial complex $\H \; = (V,H)$. We may
assume that a boolean 
representation $M$ of $\H$ is {\em
  reduced}, i.e. all the rows of $M$ are distinct and nonzero. Note
that we are assuming that $P_1(V) \subseteq H$ in all circumstances,
hence all columns must be nonzero as well.

\bl
\label{preli}
It is decidable in time $O(n!m)$ whether or not 
the set of columns of
an arbitrary $m \times n$ 
boolean matrix is independent.
\el

\proof
We use induction on $n$ to show that independence can be checked in at
most $n!m\sum_{i=0}^{n-1} \frac{1}{i!}$ elementary steps.

Assume that $n = 1$. Let $M$ denote an $m \times 1$ boolean matrix. 
Then the single column of $M$ is independent if
and only if $M$ is nonzero. Clearly, we may check if $M$ is nonzero in
$m = 1!m\sum_{i=0}^{1-1} \frac{1}{i!}$ elementary steps.

Assume now that $n > 1$ and the claim holds for $n-1$. Let $M$ denote
an $m \times n$ boolean matrix. A necessary condition for the columns
of $M$ to be independent is existence of a {\em marker} of type $j \in
\{ 1,\ldots,n\}$: a
row having a 1 at column $j$ and zeroes anywhere else. We need at most
$mn$ elementary steps to determine all $j \in \{ 1,\ldots,n\}$
admitting a marker of type $j$. For each such $j$ (and there are at
most $n$), we must check if the columns of the $(m-1)\times(n-1)$
matrix obtained by removing the marker and the $j$th column from $M$
are independent. Applying the induction hypothesis, we deduce that
independence of the columns of $M$ can be checked in at most
$$mn + n(n-1)!(m-1)\sum_{i=0}^{n-2} \frac{1}{i!} = 
\frac{n!m}{(n-1)!} + n!(m-1)\sum_{i=0}^{n-2} \frac{1}{i!} \leq 
n!m\sum_{i=0}^{n-1} \frac{1}{i!}$$
elementary steps, completing the induction. 

Since $\sum_{i=0}^{n-1} \frac{1}{i!} \leq e$, it follows that
independence can be checked on at most $en!m$ steps, hence in time $O(n!m)$.
\qed

Let $\H\; = (V,H)$ be a boolean representable simplicial complex
defined by an $R \times V$ 
boolean matrix $M = (m_{rv})$. We assume $M$ to be reduced.

For each $r \in R$, let
$$Z_r = \{ v \in V \mid m_{rv} = 0\}.$$
By \cite[Lemma 5.2.1]{RSm}, we have $Z_r \in \flatx\H$ for every $r
\in R$.

If $2 \leq |Z_r| < |V|$, then $Z_r$ is said to be a {\em line} of $M$.
We denote by $\L_M$
the set of all lines of $M$. 
Now every element of $\flatx\H$ is of 
the form $\oo{X}$ 
for some $X \in H$ by \cite[Proposition 4.2.4]{RSm}. 
On the other hand, $\oo{X} = V
\notin \L_M$ whenever $X$ is a facet of $\H$ by \cite[Proposition
4.2.4]{RSm}. It follows that 
\beq
\label{nrow}
|R| \leq |\flatx\H|-1 \leq |H \setminus \fct\H| \leq \sum_{i=0}^{d}
\binom{n}{i} \leq (d+1)n^d.
\eeq

We consider next the problem of recognizing a boolean representation
of a simplicial
complex of dimension $d \geq 0$. Note that we view $d$ as a fixed constant.

\bl
\label{recog}
Let $d \geq 0$. It is decidable in time $O(n^{2d+3})$ whether a reduced
  boolean matrix with $n$ columns defines a simplicial complex of dimension $d$.
\el

\proof
Let $M$ be such a matrix. By (\ref{nrow}), $M$ must have at most
$(d+1)n^d$ rows and we can check this necessary condition in time
$O(n^d)$, hence we may assume that $M$ has $O(n^d)$ rows. On
the other hand, $M$ has $\binom{n}{d+1}$ subsets of
$d+1$ columns. By Lemma \ref{preli}, we can decide in time
$O(n^{d})$ 
whether each such subset is a face of $\H$. Hence we can
decide in time $\binom{n}{d+1}O(n^{d})$, thus $O(n^{2d+1})$,
whether or not $\dim\H\, \geq d$. 

Since $\dim \H = d$ if and only if $\dim\H\, \geq d$ and $\dim\H\,
\not\geq d+1$, we may decide $\dim\H = d$ in time $O(n^{2d+1}) +
O(n^{2d+3})$, hence $O(n^{2d+3})$.
\qed

We present next a complexity bound for the computation of faces.

\bt
\label{detsi}
Let $d \geq 0$. It is possible to compute in time $O(n^{2d+1})$ the
list of faces of a simplicial complex of dimension $d$ defined by
a reduced boolean matrix with $n$ columns. Moreover, facets can be
marked in this list in time $O(n^{2d+2})$.
\et

\proof
Note that, by Lemma \ref{recog}, given a reduced boolean matrix $M$,
we can decide in time 
$O(n^{2d+3})$ whether $M$ defines a simplicial complex of dimension $d$. 

By (\ref{nrow}), $M$ has $O(n^d)$ rows. On
the other hand, $M$ has $\binom{n}{i}$ subsets of
$i$ columns for $i = 0,\ldots,d+1$. In view of Lemma \ref{preli}, we
can decide in time $O(n^{d})$ 
whether each such subset is a face. Hence we can
enumerate all the faces of $\H$ in time $\sum_{i=0}^{d+1}
\binom{n}{i}O(n^{d})$, thus $O(n^{2d+1})$.

For each face $I$ of dimension $< d$ and each $p \in V \setminus I$,
we can check in time $O(n^{d})$ whether $I \cup \{ p \}$ is still a
face (if $I$ has dimension $d$, is certainly a facet). Hence we may
check whether $I$ is a facet in time $O(n^{d+1})$, and so we may mark all
facets (among the $O(n^{d+1})$ faces) in time $O(n^{2d+2})$. 
\qed

We discuss now the computation of flats.

\bt
\label{gc}
Let $d \geq 2$.
It is possible to compute in time $O(n^{3d+3})$ the 
list of flats of a simplicial complex of dimension $d$ defined by
a reduced boolean matrix with $n$ columns.
\et

\proof
By Theorem \ref{detsi}, we may enumerate the list of faces
$X_1,\ldots,X_m$ of $\H$ 
in time $O(n^{2d+1})$. Note that $m \leq \sum_{i=0}^{d+1}
\binom{n}{i}$, hence $m$ is $O(n^{d+1})$.

Let $X \in H$. We claim that we can compute
$\oo{X}$ in time $O(n^{2d+3})$. Note that if $X$ is a facet, then we
have $\oo{X} = V$ by \cite[Proposition 4.2.4]{RSm}.

Indeed, let $Y = X$. By Theorem \ref{detsi}, we may check whether $Y$
contains a facet in time $O(n^{2d+2})$, yielding $\oo{Y} = V$. Hence
we may assume that $Y$ contains no facet. For every non-facet $X_i$
and $p \in V \setminus 
Y$, we may check whether $X_i \subseteq Y$ and $X_i \cup \{ p \}
\notin H$ hold simultaneously. There exist $O(n^{d})$ non-facets
$X_i$, hence we have $O(n^{d+1})$ choices for
both $i$ and $p$. Since $m$ is $O(n^{d+1})$ we may check if $X_i \cup \{ p \}
\notin H$ in time $O(n^{d+1})$. If this happens, we replace $Y$ by $Y \cup
\{ p \}$ and we restart the process. Eventually, we reach a point
where $Y$ contains a facet or there are no more $p$'s to add. In view
of \cite[Proposition 4.2.5]{RSm}, we may then deduce that $Y = \oo{X}$. 

Now each cycle $Y \mapright{} Y \cup \{ p \}$ can be performed in time
$O(n^{2d+2})$ and there are at most $n$ cycles to be performed, hence
$\oo{X}$ can be computed in 
time $O(n^{2d+3})$. Since the number of non-facets $X_i$ is
$O(n^{d})$, we can compute their closures (and consequently all flats)
in time $O(n^{3d+3})$.  
\qed



\bc
\label{cgf}
Let $d \geq 2$. Let $\H$ denote an
arbitrary simplicial complex of 
dimension $d$ represented by a reduced boolean matrix $M$ with $n$
columns. Then: 
\bi
\item[(i)] $\Gamma{\rm Fl}\H$ can be
  computed in time 
$O(n^{2d+5})$;
\item[(ii)] $\pi_1(||\H||)$ can be computed in time
  $O(n^{2d+5})$.
\ei
\ec

\proof
(i) We have $\binom{n}{2}$ potential edges $a \edge b$ in
$\Gamma\flatx\H$. By the proof of Theorem \ref{gc}, we may compute
$\oo{ab}$ in time $O(n^{2d+3})$, and check whether or not $\oo{ab} = V$.
Thus we reach a global complexity bound of $O(n^{2d+5})$.

(ii) By Theorem \ref{fung}, we need to compute the number of
connected components of $\Gamma\flatx\H$ (a graph with $n$ vertices
and at most $\binom{n}{2}$ edges) and to identify the $H$-trivial 
components. It is easy to see by induction that the number of
connected components can be computed in time $O(n^2)$. In view of
Theorem \ref{detsi}, we can identify the $H$-trivial connected
components in time $O(n^{2d+3})$. Therefore
$\pi_1(||\H||)$ can be computed in time $O(n^{2d+5}) + O(n^2) +
O(n^{2d+3}) = O(n^{2d+5})$.
\qed

We show next how these complexity bounds can be improved in the case
of dimension 2. 

Let $\Gamma = (V,E)$ be a graph. Given $v \in V$, we write ${\rm
  nbh}(v) = \{ w \in V \mid vw \in E \}$. We say that $A \subseteq V$
is a {\em superanticlique} if $|A| > 1$ and
$$\nbh(a) \cup \nbh(b) = V \setminus A$$
holds for all $a,b \in A$ distinct.
In particular, the superanticlique $A$
is a maximal anticlique (i.e. maximal with respect to $P_2(A) \cap E =
\emptyset$).

Superanticliques play a major role in the theory of
boolean representable simple simplicial complexes of dimension 2. Let
$M$ be a boolean matrix representation of such a
complex, say $\H\; = (V,H)$. We denote by $\Gamma M$ the graph with
vertex set $V$ and edges of the form $p \edge q$ whenever $pq$ is a
2-subset of a line of $M$. By
\cite[Theorem 6.3.6]{RSm},  
$\flatx\H$ is the union of $P_{\leq 1}(V) \cup \{ V \} \cup \L_M$ with
the set of all superanticliques of $\Gamma M$.

Given two graphs $\Gamma = (V,E)$ and $\Gamma' = (V',E')$, assumed to
be disjoint, we define their {\em join} to be the graph $\Gamma +
\Gamma' = (V \cup V',E \cup E' \cup E'')$, where $E'' = \{ vv' \mid v
\in V, v' \in V' \}$. Their {\em coproduct} is the graph $\Gamma
\sqcup \Gamma' = (V \cup V',E \cup E')$.

Given $n \geq 1$, we denote by $K_n$ the complete graph on $n$
vertices. We denote by $\oo{K_n}$ the complement graph of $K_n$, so
that $\oo{K_n}$ has $n$ vertices and no edges.

We define now two classes of graphs as follows. Let $\Omega_1$ be the
class of all graphs of the form $(\oo{K_n} + \Delta) \sqcup K_1$,
where $n \geq 1$ and $\Delta$ is any finite graph. 
Let $\Omega_2$ be the
class of all graphs of the form $(K_1 + \Delta) \sqcup (K_1 + \Delta')$,
where $\Delta$ and $\Delta'$ are any finite graphs. 

\bt
\label{mtoh}
Let $M$ be a boolean matrix representation of a simple simplicial
complex $\H$ of dimension 2. Then:
\bi
\item[(i)] if $\Gamma M$ is connected or belongs to
  $\Omega_1 \cup \Omega_2$, then $\Gamma{\rm Fl}\H$ is connected;
\item[(ii)] in all other cases, $\Gamma{\rm Fl}\H = \Gamma M$.
\ei
\et

\proof
(i) Since $\L_M \subseteq \flatx\H$ by \cite[Lemma
  5.2.1]{RSm}, then
$\Gamma M$ is a subgraph of $\Gamma\flatx\H$ with the same
vertex set. Therefore $\Gamma M$ connected implies
$\Gamma\flatx\H$ connected.

Assume next that $\Gamma M \in \Omega_1$, say of the form
$(\oo{K_n} + \Delta) \sqcup K_1$. Let $A$ be the union of the $n$
vertices of $\oo{K_n}$ and the single vertex of $K_1$. Given $a,b \in
A$, then $\nbh(a) \cup \nbh(b)$ are the vertices of $\Delta$, i.e. $V
\setminus A$. Thus $A$ is a superanticlique of $\Gamma M$ and so $A
\in \flatx\H \setminus \{ V \}$ by \cite[Theorem 6.3.6]{RSm}. Since
$A$ intersects the two connected components of $\Gamma M$, it follows
that $\Gamma\flatx\H$ is connected.

Assume now that $\Gamma M \in \Omega_2$, say of the form
$(K_1 + \Delta) \sqcup (K_1 + \Delta')$.
Let $A$ consists of the two 
vertices in both copies of $K_1$, say $a,b$.
Then $\nbh(a) \cup \nbh(b)$ are the vertices of $\Delta$ and $\Delta'$, i.e. $V
\setminus A$. Thus $A$ is a superanticlique of $\Gamma M$ and so $A
\in \flatx\H \setminus \{ V \}$. Since
$A$ intersects the two connected components of $\Gamma M$, it follows
that $\Gamma\flatx\H$ is connected.

(ii) Suppose that $\Gamma M$ is disconnected and $\Gamma{\rm Fl}\H
\neq \Gamma M$. We must show that $\Gamma M \in \Omega_1 \cup \Omega_2$.
In view of \cite[Theorem 6.3.6]{RSm}, there exists
some superanticlique $A$ of $\Gamma M$. It follows from the definition
that $A$ must intersect all the connected components of $\Gamma M$.

Suppose that $\Gamma M$ has
more than two connected components. Since $\H$ has dimension 2, one of
the connected components, say $C$, must be nontrivial. Let $a,b \in A
\setminus C$. Then $(\nbh(a) \cup \nbh(b))\cap C = \emptyset$. 
Since $C \setminus A \neq \emptyset$, this contradicts $\nbh(a) \cup
\nbh(b) = V \setminus A$. Therefore $\Gamma M$ has precisely two
connected components, and we may write $\Gamma M = \Gamma \sqcup
\Gamma'$ with $\Gamma$ and $\Gamma'$ connected. 

Suppose that $\Gamma$ and $\Gamma'$ are both nontrivial. The same
argument used above implies that $A$ has one element $a$ in $\Gamma$ and
another $b$ in $\Gamma'$. Since $\nbh(a) \cup \nbh(b) = V \setminus \{
a,b\}$, it follows that $\Gamma M \in \Omega_2$.

Thus we may assume that $\Gamma'$ is trivial. Let $a \in A$ be a
vertex of $\Gamma$ and let $b$ be the unique vertex of $\Gamma'$
(which is in $A$). Let $\Delta$ (respectively $\Delta'$) be the
subgraph of $\Gamma$ induced by 
$\nbh(a)$ (respectively the remaining vertices of $\Gamma$). Since $A
= V \setminus (\nbh(a) \cup \nbh(b))$, then $\Delta'$ is an edgeless
graph. Let $c$ be a vertex of $\Delta'$. Since $\nbh(c) \cup \nbh(b) =
V \setminus A = \nbh(a) \cup \nbh(b)$, it follows that $\Gamma =
\Delta + \Delta'$. Therefore $\Gamma M \in \Omega_1$.
\qed

Now we can provide complexity bounds for both fundamental group and
decidability of shellability in dimension 2.  

\bt
\label{ccb}
Let $\H$ denote an arbitrary simplicial complex of 
dimension $2$ represented by a reduced  boolean matrix $M$ with $n$
columns. Then: 
\bi
\item[(i)] if $\H$ is simple, then $\pi_1(||\H||)$ can be
  computed in time
  $O(n^4)$;
\item[(ii)] it can be determined in time $O(n^4)$ whether or not $\H$
  is shellable. 
\ei
\et

\proof
(i) Since $\H$ is connected by Lemma \ref{connected}, then
$\pi_1(||\H||)$ is well defined. 
Since $M$ is reduced, it has at most $3n^2$ rows by (\ref{nrow}).

By Corollary \ref{simfun}, it suffices to compute the
number of connected components of $\Gamma\flatx\H$. 

Since $M$ has at most $3n^2$ rows, we may compute $\Gamma M$ in time
$O(n^4)$ (there are $\binom{n}{2}$ pairs of vertices to
check, and each pair can be checked in time $O(n^2)$). 

We claim that we can check whether or not $\Gamma M$ is connected in
time $O(n^4)$. Indeed, let $r$ be the number of rows of $M$ and let
$M_i$ be the submatrix of $M$ defined by the first $i$ rows $(i =
1,\ldots,r)$. Obviously, we can compute the connected components of
$\Gamma M_1$ in time $O(n)$. Assume now that $1 < i \leq r$ and the 
connected components of $\Gamma M_{i-1}$ 
were computed in time $O(in^2)$. We can mark the zero
entries of the $i$th row with the connected components of 
$\Gamma M_{i-1}$ in time $O(n^2)$ and merge distinct connected
components arising this way in time $O(n^2)$, and the complexity
constants for these two procedures do not depend on $i$. Since $r \leq
3n^2$, it follows by induction that the connected components of
$\Gamma M_{r} = \Gamma M$ 
can be computed in time $O(n^4)$. Therefore we can check whether or
not $\Gamma M$ is connected in time $O(n^4)$.

We claim that we can also decide whether or not $\Gamma M \in \Omega_1
\cup \Omega_2$ in time $O(n^4)$. Since the connected components of 
$\Gamma M$ were already computed in time $O(n^4)$, it suffices to show
that it is decidable in time $O(n^4)$ whether or not a connected graph
with at most $n$ vertices is of the form $K_1 + \Delta$ or $\oo{K_m} +
\Delta$. The first case is obvious since we have at most $n$ potential
choices for the vertex playing the $K_1$ role. For the case $\oo{K_m}
+ \Delta$, we note that we need at most $n$ tries to pick a vertex $v$
in $\oo{K_m}$, and for each such $v$ the vertices of $\Delta$ (if it
exists) would be necessarily $\nbh(v)$, hence the vertices in both
$\oo{K_m}$ and $\Delta$ would be fully determined by $v$. We would be
able to mark them as such in time $O(n)$. Finally, we may decide
whether $\nbh(v)$ is an anticlique in time
$O(n^2)$, and we can check whether $a \edge b$ is an edge for all $a
\in \nbh(v)$ and $b \notin \nbh(v) \cup \{ v \}$ in time $O(n^2)$,
proving our claim.

Now it follows from Theorem \ref{mtoh} that we may compute the number
of connected components of $\Gamma\flatx\H$ in time $O(n^4)$, and we
apply Theorem \ref{fung}.

(ii) By Proposition \ref{scon}, we can produce a submatrix $M'$
of $M$ representing $\H_S$ by removing repeated columns. We may do it
by comparing pairs of columns. There are $\binom{n}{2}$ pairs to
compare, and each pair can be compared in time $O(n^2)$, hence we can
compute $M'$ in time $O(n^4)$.

In view of Theorem \ref{shelltwo}, we can assume that
$\H$ is simple, and use the proof of part (i).
\qed

Note that the quartic bound in part (i) is much better than the
$O(n^9)$ bound provided by Corollary \ref{cgf}(i).

We remark also that, once shellability is ensured, an actual
shelling can be produced in the simple case using the algorithms
described in \cite[Lemma 7.2.1]{RSm} and \cite[Lemma 7.2.5]{RSm} within the
same quartic complexity bounds. The
extension to the general case follows then from Proposition
\ref{retro}(x) and Theorem \ref{shelltwo}. Therefore we obtain the
following corollary.

\bc
\label{consh}
Let $\H$ denote an arbitrary shellable simplicial complex of 
dimension $2$ represented by a reduced  boolean matrix $M$ with $n$
columns. Then a shelling of $\H$ can be actually computed in time
$O(n^4)$. 
\ec

The $i$-th Betti number $w_i(\H)$ is defined as the rank of the $i$th
homology group of $||\H||$. If $\H$ is shellable, then
by \cite{BW} $w_i(H)$ is the number of
homology facets in a shelling $B_1,\ldots,B_t$ of $\H$. 
We say that $B_k$ $(k > 1)$ is a {\em homology facet} in this 
shelling if $2^{B_k} \setminus \{ B_k \} \subseteq \cup_{i = 1}^{k-1}
2^{B_i}$.

Assume that $\H$ satisfies the conditions of Corollary
\ref{consh}. Then we can construct a shelling $B_1,\ldots,B_t$
in time $O(n^4)$. Now we can build a sequence
$\Delta_1,\ldots,\Delta_t$ of graphs with vertex set $V(\Delta_k) =
\cup_{i=1}^k B_i$ and edge set
$E(\Delta_k) = \cup_{i=1}^k P_2(B_i)$ to help us keep track of
homology facets: indeed, if $k > 1$, then $B_k$ is a homology facet if
and only if ($|B_k| = 2$ and $B_k \subseteq V(\Delta_{k-1})$) or
($|B_k| = 3$ and $P_2(B_k) \subseteq E(\Delta_{k-1})$). 
Since $t \in O(n^3)$, this provides a proof for the following result.


\bc
\label{betti}
Let $\H$ denote an arbitrary shellable simplicial complex of 
dimension $2$ represented by a reduced  boolean matrix $M$ with $n$
columns. Then the Betti numbers of $\H$ can be computed in time
$O(n^4)$. 
\ec

\section{Open problems}

The problem of determining the homotopy type for BRSCs of dimension
$\geq 3$ remains open, as are the problems of identifying the
shellable and the sequentially Cohen-Macauley BRSCs for dimension
$\geq 3$.

\section*{Acknowledgments}

The first author was partially supported by Binational Science
Foundation grant number 2012080.

The second author thanks the Simons Foundation-Collaboration Grants
for Mathematicians for travel grant $\#$313548.

The third author was partially supported by
CNPq (Brazil) through a BJT-A grant (process\linebreak
313768/2013-7) and
CMUP (UID/MAT/00144/2013),
which is funded by FCT (Portugal) with national (MEC) and European
structural funds through the programs FEDER, under the partnership
agreement PT2020.

\bigskip

{\sc Stuart Margolis, Department of Mathematics, Bar Ilan University,
  52900 Ramat Gan, Israel} 

{\em E-mail address:} margolis@math.biu.ac.il

\bigskip

{\sc John Rhodes, Department of Mathematics, University of California,
  Berkeley, California 94720, U.S.A.}

{\em E-mail addresses}: rhodes@math.berkeley.edu, BlvdBastille@aol.com

\bigskip

{\sc Pedro V. Silva, Centro de
Matem\'{a}tica, Faculdade de Ci\^{e}ncias, Universidade do
Porto, R. Campo Alegre 687, 4169-007 Porto, Portugal}

{\em E-mail address}: pvsilva@fc.up.pt


\begin{thebibliography}{99}
\bibitem{Bjo2} A.~Bj\"orner, Shellable and Cohen-Macaulay partially
  ordered sets, {\em Trans. Amer. Math. Soc.} 260 (1980), 159--183.
\bibitem{Bjo} A.~Bj\"orner, The homology and shellability of matroids
  and geometric lattices, In: {\em Matroid applications},
  Encycl. Math. Appl. 40, Cambridge University Press, 1992, pp. 226--283.
\bibitem{BW} A.~Bj\"orner and M.~L.~Wachs, Nonpure shellable complexes
  and posets I, {\em Trans. Amer. Math. Soc.} 348 (1996), 1299--1327.
\bibitem{BW2} A.~Bj\"orner and M.~L.~Wachs, Nonpure shellable complexes
  and posets II, {\em Trans. Amer. Math. Soc.} 349 (1997), 3945--3975.
\bibitem{BWW} A.~Bj\"orner, M.~Wachs and V.~Welker, On sequentially
  Cohen-Macaulay complexes and posets, {\em Israel J. Math.} 169
  (2009), 295--316.
\bibitem{Duv} A.~M.~Duval, Algebraic shifting and sequentially
  Cohen-Macaulay simplicial complexes, {\em Electron. J. Combin.} 3
  (1996).
\bibitem{Hat} A.~Hatcher, {\em Algebraic Topology}, Cambridge
  University Press, Cambridge, 2002.
\bibitem{Izh2}
Z.~Izhakian, Tropical arithmetic and tropical matrix
algebra, {\em Comm. Algebra} 37(4) (2009), 1--24.
\bibitem{IR1}
 Z.~Izhakian and J.~Rhodes, New representations of matroids and
 generalizations, preprint, arXiv:1103.0503, 2011.
\bibitem{IR2}
 Z.~Izhakian and J.~Rhodes, Boolean representations of matroids and
 lattices, preprint, arXiv:1108.1473, 2011.
\bibitem{IR3}
 Z.~Izhakian and J.~Rhodes, C-independence and c-rank of posets and
 lattices, preprint, arXiv:1110.3553, 2011.
\bibitem{IRow2}
 Z.~Izhakian and L.~Rowen, The tropical rank of a tropical matrix, {\em
 Commun. Algebra} 37(11) (2009), 3912--3927.
\bibitem{IRow}
 Z.~Izhakian and L.~Rowen, Supertropical algebra, {\em Adv. Math.}
 225(8) (2010), 2222--2286.
\bibitem{IRow22}
 Z.~Izhakian and L.~Rowen, Supertropical matrix algebra, {\em Israel
 J. Math.} 182 (2011), 383--424.
\bibitem{LS} R.~C.~Lyndon and P.~E.~Schupp, {\em Combinatorial Group 
Theory}, Springer-Verlag, 1977.
\bibitem{Oxl} J.~G.~Oxley, {\em Matroid Theory}, Oxford Science
   Publications, 1992.
\bibitem{RSm} J.~Rhodes and P.~V.~Silva, {\em Boolean Representations
  of Simplicial Complexes and Matroids}, Springer Monographs in
  Mathematics, Springer, 2015.
\bibitem{Rot} J.~J.~Rotman, {\em An Introduction to Algebraic
  Topology}, Springer-Verlag, 1986.
\bibitem{Spa} E.~Spanier, {\em Algebraic Topology}, Corrected reprint,
  Springer-Verlag, New York-Berlin, 1981.
\bibitem{Sta2} R.~P.~Stanley, {\em Combinatorics and Commutative
    Algebra}, 2nd edition, Birkh\"auser, Boston, 1995.
\bibitem{Wac} M.~Wachs, Poset topology: tools and applications, {\em
   Geometric Combinatorics} 13 (2007), 497--615.
\bibitem{Wal} C.~T.~C.~Wall, Finiteness conditions for CW-complexes,
  {\em Ann. Math.} 81(1) (1965), 56--69.
\bibitem{Whi2} N.~White (ed.) {\em Theory of Matroids}, Encyclopedia of
  Mathematics and its Applications 26, Cambridge University Press, 1986.
\end{thebibliography}
\end{document}